\documentclass[11pt,leqno]{article}
\usepackage{lipsum}
\usepackage{amsfonts}
\usepackage{graphicx}
\usepackage{epstopdf}
\usepackage{algorithmic}
\usepackage{dsfont}
\usepackage{amssymb}

\usepackage[T1]{fontenc}
\usepackage{amsmath}
\usepackage{amsthm}
\usepackage{color}
\usepackage{graphicx}
\usepackage{lmodern}
\usepackage{mathrsfs}
\usepackage{appendix}

\usepackage{url}
\ifpdf
  \DeclareGraphicsExtensions{.eps,.pdf,.png,.jpg}
\else
  \DeclareGraphicsExtensions{.eps}
\fi

\newcommand{\cG}{{\mathcal G}}
\newcommand{\cP}{{\mathcal P}}
\newcommand{\cPmzeroGammaLipCV}{\cP_{m_0}(\Gamma^{\textrm{Lip}}_{C,V})}
\newcommand{\GammamuoptC}  {\Gamma^{\mu,{\rm opt}}_C}
\newcommand{\R}{{\mathbb R}}

\newcommand{\N}{{\mathbb N}}
\newcommand{\cD}{{\mathcal D}}
\newcommand{\cL}{{\mathcal L}}
\let\G=\Gamma
\let\d=\delta

\newcommand{\vertex}{{0}}

\newcommand{\car}{{\mathds{1}}}
\newcommand{\conv}{{\mathrm{conv}}}
\let\ds=\displaystyle

\numberwithin{equation}{section}

\newtheorem{theorem}{Theorem}[section]
\newtheorem{lemma}[theorem]{Lemma}
\newtheorem{proposition}[theorem]{Proposition}
\newtheorem{definition}{Definition}[section]

\newtheorem{remark}[theorem]{Remark}

\title{Deterministic mean field games on networks\\ and \\ related optimal control problems  
  }

\author{Yves Achdou\thanks{Universit\'e de Paris Cit\'e and Sorbonne Universit\'e, CNRS, Laboratoire Jacques-Louis Lions, (LJLL), F-75006 Paris, France, achdou@ljll-univ-paris-diderot.fr
    .}
  \and Claudio Marchi \thanks{Dipartimento di Matematica ``Tullio Levi-Civita'', Universit\`a di Padova, 35121 Padova, Italy, marchi@math.unipd.it
    .}
\and Nicoletta Tchou\thanks{Univ Rennes, CNRS, IRMAR - UMR 6625, F-35000 Rennes, France, nicoletta.tchou@univ-rennes.fr
  .}}

\usepackage{amsopn}

\begin{document}

\maketitle

\begin{abstract}
We study a class of deterministic mean field games and related optimal control problems, with a finite time horizon and in which the state space is a network.

An agent controls her velocity, and, when she occupies a vertex, she can either remain still or enter any adjacent edge. The running and terminal costs are assumed to be continuous in each edge, but
may jump at the vertices.
%not necessarily globally continuous on the network, i.e., 
%the costs  
Compared to the companion paper \cite{AMMT3}, we make more general assumptions about the costs and consider networks with an arbitrary number of vertices; this higher degree of generality brings new difficulties.

For the optimal control problems mentioned above, we obtain in particular the existence of optimal trajectories and regularity results concerning the optimal trajectories and the value function.

These control theoretic results make it possible
to address a class of mean field games on networks, with costs that do not depend separately on the control and on the distribution of states,
and that are non-local with respect to the latter. Focusing on a Lagrangian formulation, we obtain the existence of relaxed equilibria consisting of probability measures on admissible trajectories. To any relaxed equilibrium corresponds a mild solution, i.e. a pair $(u, m)$ made of the value function $u$ of a related optimal control problem and a family $m = (m(t))_t$ of probability measures on the network. Given $m$, the value function $u$ is a viscosity solution of a Hamilton-Jacobi problem on the network. We then investigate the regularity properties of $u$ and a weak form of a Fokker-Planck equation satisfied by $m$.

\end{abstract}


\paragraph{MSCcodes}
35F50, 35Q89, 35Q91, 35R02, 49K20, 49L25, 49N80, 91A16.

\section{Introduction}\label{intro}
The present paper is devoted to deterministic Mean Field Games (MFGs in short)  and related optimal control problems, in which the state space is a network, i.e. a connected set obtained as the union of a finite number of lines, the intersections of which are named vertices.
 Given the time evolution of the distribution of agents in the state space, each agent solves an optimal control problem with a finite time horizon, the control variable being her velocity. In particular, when an agent is at a vertex, she can either remain still or enter any adjacent edge. The running and terminal costs depend on the distribution of states in a non-local, regularizing manner but are not supposed to be continuous across the vertices (the definition of the cost functions may change from one edge to the next). The present paper is a companion paper to \cite{AMMT3}, in which P. Mannucci and the authors addressed similar MFGs, but focused on the situation in which the network has only one vertex, the running cost depends separately on the control variable and on the distribution of states and is purely quadratic in the control variable.
Here, compared to \cite{AMMT3}, the assumptions made on the running costs are more general, and general networks are investigated in Section \ref{sec_general_nets} below. In particular, the costs are no longer assumed to depend separately on the control variable and on the distribution of states.
To keep the present paper reasonably short, we shall avoid
repeating the arguments of \cite{AMMT3} as much as possible and focus on the new arguments. 

In the most commonly studied situations, a deterministic mean field game is described by the pair made of the distribution of states at all times and the optimal value of a representative agent, and the latter quantities satisfy boundary value problem involving a system of PDEs coupling a continuity equation (forward in time) and a Hamilton--Jacobi (HJ) equation (backward in time), see the pioneering articles of Lasry and Lions \cite{LL1,LL2,LL3} and the notes written by Cardaliaguet
\cite{C}.
%and Porretta in \cite{MR4214773}. 
However, there are cases in which it is difficult to characterize the MFG by a boundary value problem.
For example, it was observed in \cite{AHLLM} that if state constraints are imposed, the distribution of states may become singular near the boundary of the state space, and it becomes difficult to describe its evolution by a PDE
supplemented with boundary conditions.
For this reason, in order to study MFGs in which the agents are constrained to remain in the closure of a regular domain of $\R^N$, Cannarsa et al., following ideas contained in \cite{Br89,Br93,BB00,BC15,CM16}, introduced a notion of relaxed equilibrium defined in a Lagrangian setting rather than with PDEs, see \cite{CC,CCC1,CCC2}. The evolution of the game is described by probability measures defined on a set of admissible trajectories rather than by time-dependent probability measures defined on the state space. Similar strategies have been followed by Mazanti and Santambrogio in \cite{MS19} to study MFGs in which every agent aims at exiting a given closed subset of a general compact metric space in minimal time, and by P. Mannucci with the authors in \cite{AMMT2},  to study deterministic state constrained MFGs in which the agents control their acceleration (a non-locally controllable case).

In \cite{AMMT3} and here,  we also deal with relaxed equilibria defined in a Lagrangian setting rather than with PDEs. The main new difficulties compared to \cite{CC,CCC1,CCC2} come from the geometry which is intrinsically singular, and also from the discontinuities of the costs at the vertices. In \cite{AMMT3}, these features explain why singularities in the state distribution may appear at the vertices, then propagate into the edges of the network, and possibly split when crossing a vertex or simply disappear.  The present situation may be even worse because it cannot be ruled out that singularities develop inside edges.
However, with each relaxed equilibrium,  it will be possible to  associate a pair $(u,m)$, where $u$ is a viscosity solution to a Hamilton-Jacobi equation posed in the network, while $m$ is a solution to a continuity equation posed on the network, in a weak sense that will  be defined.
We think that the assumptions made in the present paper are close to minimal to obtain the existence of optimal trajectories and the boundedness of optimal controls. However, because of their generality, they will not allow us to recover the most accurate results contained in \cite{AMMT3} on the propagation of singularities in the state distribution, see Remarks \ref{rmk:h7} and \ref{rem:nodirac} below.

Our main strategy is as follows. First, we establish control-theoretic results that will be important in the further analysis of the mean field game.  Next, we investigate relaxed mean field equilibria relying on the latter results.
Then, given a relaxed equilibrium, we analyze the related pair (value function, flow of probability measures) on the network, applying the previously obtained regularity results to the related optimal control problem. Compared to \cite{AMMT3}, the higher degree of generality brings new difficulties, for example
when proving the existence of optimal trajectories (this is important for defining relaxed equilibria) or when obtaining regularity properties of the value function (which will later be used for studying the evolution of the distribution of the agents in the state space). More precisely, in  \cite{AMMT3}, the existence of optimal trajectories mostly stems from the lower semi-continuity of the cost with respect to the $L^2$-weak convergence of controls (that is obtained in a straightforward manner); here, the  lower semi-continuity of the cost is not obvious and more involved arguments are needed: we extend the cost to a functional defined on trajectories in a finite dimensional vector space and prove several properties of this new cost (measurability, convexity, etc...), that finally allow us to invoke a classical result in the calculus of variations.  Concerning the regularity of optimal trajectories, we will see that the Euler-Lagrange first-order optimality conditions comprise ordinary differential equations in the interior of the edges; in \cite{AMMT3}, the latter equations are in normal form and the assumptions of Cauchy-Lipschitz theorem hold, but this is not the case here: the differential equations are not in normal form and the uniqueness of local solutions is not guaranteed. This issue explains why we need to completely change the proof of the Lipschitz continuity of optimal trajectories, and carefully revisit the arguments regarding the local semi-concavity of the value function, optimal synthesis, and the continuity equation arising in mean field games.

The present paper is organized as follows.  The remainder of Section \ref{intro} contains the description of the geometry and the definition of some notation.  Section~\ref{sec:OC} is devoted to optimal control problems on a single-vertex network (which arise if the distribution of states in the MFG is given)
 and contains the following results:
the existence of an optimal trajectory for a given initial state, a closed graph property for the map which associates to each point on the network the set of optimal trajectories starting from that point, necessary optimality conditions for an optimal control, the boundedness of the optimal controls and the local or global Lipschitz regularity of the value function.  We will also see that the value function is a generalized viscosity solution of a Hamilton-Jacobi problem posed on the network with suitable conditions at the vertex.
Section~\ref{sect:MFGequil} deals with relaxed equilibria for MFGs on the same single-vertex network and relies on the results of Section~\ref{sec:OC}. As in \cite{AMMT3}, the existence of a relaxed equilibrium is obtained using Kakutani's fixed-point theorem and the closed graph property mentioned above.  Then it is possible to associate to a relaxed equilibrium a family of time-dependent probability measures on the state space $(m(t))_t$ and the value function $u$ of an optimal control problem involving $m$.  The pair $(u,m)$ is named a {\sl mild solution} of the MFG, see \cite{CC}. All the results in  Section~\ref{sec:OC}  apply to the latter optimal control problem. In particular, we deduce some regularity properties of $u$ and prove that $m$ solves a continuity equation in a weak sense. Finally, in Section \ref{sec_general_nets}, we show how to deal with a general network (with several vertices).
\subsection{Notation}\label{sec:notations}

\paragraph{Junctions}
A {\sl junction} $\cG$ is the union of $N>1$ semi-infinite lines $J_i$ with a common endpoint $\vertex$.
The lines $J_i$ will be referred to as {\sl the edges}, while the point $\vertex$ will be referred to as {\sl the vertex}. Because only metric properties will matter for the optimal control problem, it is not restrictive to see the junction $\cG$ as a subset of $\R^N$, set the origin of $\R^N$ at $\vertex$, and assume that the edges $J_i$ are the closed half-lines $J_i=\R_+ e_i$, where $(e_i)_{1\le i\le N}$ is the canonical basis of $\R^N$.
The geodesic distance $d(x,y)$ between two points $x,y$ of $\cG$ is 
\begin{displaymath}
  d(x,y)=\left\{\begin{array}[c]{ll}
|x-y|\quad \hbox{if } x,y \hbox{ belong to the same edge } J_i \\
|x|+|y|  \quad \hbox{if } x,y \hbox{ belong to different edges } J_i \hbox{ and } J_j.
\end{array}\right.
\end{displaymath}
              
\paragraph{The set  $C^1(\cG)$.} Let  $C^1(\cG)$ be  the set of functions $\varphi\in C(\cG)$ such that, 
for every $i=1,\dots,N$, the restriction of $\varphi$ to the edge~$J_i$, $\varphi_{\mid J_i}$ belongs to $C^1(J_i)$. For $\varphi\in C^1(\cG)$,  set
\begin{equation}\label{eq:def_deriv}
D\varphi(x)=\left\{\begin{array}{ll}
\partial_{x_i} \varphi_{\mid J_i}(\bar x)  &\qquad\textrm{if }x=\bar x e_i \in J_i\setminus\{\vertex\},\\
\left(\partial_{x_1} \varphi_{\mid J_1}(0),\dots,\partial_{x_N} \varphi_{\mid J_N}(0)
\right)&\qquad\textrm{if }x=\vertex.
\end{array}\right.
\end{equation}
Observe that $D\varphi(x)$ is a scalar when the point $x$ lies in the interior of a given edge while it is $N$-dimensional when $x$ coincides with the vertex~$\vertex$.
When $\varphi$ also depends on the time variable, we will write $D_x$ instead of $D$.

Similarly, let $C^1(\cG\times [0,T])$ be the set of functions $\varphi\in C(\cG\times [0,T])$ such that for every $1\le i\le N$, the restriction $\varphi_{\mid J_i\times [0,T]}$ belongs to $C^1(J_i\times[0,T])$.

\paragraph{Other notation}

For a metric space $A$,  $C_b(A)$ denotes the space of bounded and continuous real-valued functions on~$A$.

For any function $f: \R^d\to (-\infty,+\infty]$, $\conv f$ stands for the {\it convex envelope} of $f$, namely the largest convex function bounded from above by $f$ (see \cite[page 36]{MR1451876}). Moreover, for any collection of functions $\{f_i\}_{i\in I}$ defined in $\R^d$ with values in $(-\infty,+\infty]$,  the convex envelope $\conv \{f_i\}_{i\in I}$ is the largest convex function bounded from above by every~$f_i$. Recall that $\conv \{f_i\}_{i\in I}=\conv \left(\inf_{i\in I}f_i\right)$.\\

\section{A class of optimal control problems on $\cG$}\label{sec:OC}

In Section \ref{sect:MFGequil} below, we are going to focus on a class of MFGs on $\cG$ in which the agents
control their velocity and pay a running cost which depends locally on the state and control variables and non-locally on
the distribution of states. Here, we address the optimal control problems obtained from the latter MFGs by supposing that the distribution of states is given at all times. The results found for this class of optimal control problems will be used in Section \ref{sect:MFGequil}.
\subsection{Setting and assumptions} \label{setting}
In the remainder of Section \ref{sec:OC}, we will always assume that
the dynamics and costs satisfy the hypotheses made in the following paragraphs and will not repeat them. We will specify when additional hypotheses are needed.
\paragraph{Admissible trajectories}
Fix $T>0$. Consider a real number $p>1$.
Given $x\in \cG$ and a time $t$, $0\le t\le T$, the admissible trajectories starting from $x$ at time $t$ are functions
$y\in C([t, T]; \cG)\cap W^{1,p} (t, T; \R^N)$ such that $y(t)=x$. Given an admissible trajectory $y$, note that the sets $\Bigl\{s\in [t,T]: y(s)\in J_i\setminus\{\vertex\}\Bigr \}$ are relatively open in $[t,T]$ and that the set $\{s\in [t,T]: y(s)=\vertex\}$ is closed, so they all are measurable. Let $y'$ be the time derivative of $y$ (defined in the weak sense) on $(t, T)$. 
%The notation $\alpha$ for $y'$ is intended to stress the fact that the control variable is precisely the velocity. 
%The identity 
%$y(s)=   x+ \int_t^s y'(\theta) d\theta$ holds pointwise.
By Stampacchia's theorem (see, for instance, \cite[Lemma 7.7]{GT}), $y'(s) =0$ at almost all $s$ such that $y(s)=\vertex$. On the other hand, $y'(s)$ is aligned with $e_i$ at almost all $s$ such that $y(s)\in J_i\setminus\{\vertex\}$. Hence, it is possible to define in a unique way $\bar y'\in L^p(t,T) $ by
\begin{eqnarray}
\label{eq:baralpha1}
\bar y'(s) =0, \quad \hbox{for a.a. } s\in [t,T] \hbox{ s.t. } y(s)=\vertex,\\
\label{eq:baralpha2}
y'(s)= \bar y'(s) \sum_{i=1}^N    \car_{y(s)\in J_i\setminus\{\vertex\}}  e_i,\quad \hbox{for a.a. } s\in [t,T].
\end{eqnarray}
Given $t\in [0,T]$, let $\Gamma_t[x]$ be the set of admissible trajectories starting from $x$ at time $t$. For $t_1, t_2$,  $0\le t_1\le t_2\le T$, let $\Gamma_{t_1,t_2}[x]$ be the set of restrictions of the functions in $\Gamma_{t_1}[x]$ to the interval $[t_1,t_2]$. To alleviate the notation, we drop the subscript $t$ in $\Gamma_t[x]$ when $t=0$, i.e., we set $\Gamma[x]=\Gamma_0[x]$.

Let $\Gamma$ be the set of all admissible trajectories on $[0,T]$, i.e., $ \Gamma=\cup_{x\in \cG} \Gamma[x]$.

For any $y\in \Gamma_t[x]$, we will often set $\alpha=y'\in L^p(t,T;\R^N)$, $\bar y =|y|\in W^{1,p}(t,T;\R_+)$ and $\bar\alpha= \bar y '\in L^p (t,T)$, for which \eqref{eq:baralpha1} and \eqref{eq:baralpha2} hold. We will use the same notation for $y\in \Gamma_{t_1,t_2}[x]$ with obvious modifications in the definitions.

We skip the proof of the following lemma.
\begin{lemma}\label{compactness_lemma}
  The set $\Gamma_t[x]$ is closed with respect to the weak topology of $W^{1,p}( t,T; \R^N)$.
  The bounded subsets of $\Gamma_t[x]$ (in $W^{1,p} (t, T; \R^N)$) are relatively compact with respect to the weak topology of $W^{1,p}( t,T; \R^N)$ and in $C( [t,T]; \cG)$.
\end{lemma}

\begin{remark}[Concatenation of two admissible trajectories]\label{rmk:concat}
For $0\leq t_1\leq t_2\leq t_3\leq T$ and $x\in\cG$, if $y_1\in\Gamma_{t_1,t_2}[x]$ and $y_2\in\Gamma_{t_2,t_3}[y_1(t_2)]$,  then $\tilde y$ defined by $\tilde y(s)=y_1(s) \car_{s\in[t_1,t_2]} +y_2(s) \car_{s\in(t_2,t_3]}$ belongs to~$\Gamma_{t_1,t_3}[x]$.
\end{remark}

\paragraph{The costs}

Given a number $g_*$ and $N$ functions $g_i\in C_b(\R_+)$, the terminal cost is defined on $\cG$ by
\begin{equation}
 \label{eq:461}
g(y)= \sum_{i=1}^N g_i(|y|)\car_{y\in J_i\setminus\{\vertex\}}+\min\left\{g_*,\min\limits_{i=1,\dots,N}g_i(0)\right\}\car_{y=\vertex}.  
\end{equation}

For $i=1,\dots,N$, let us introduce the running cost $\ell_i$ associated with the edge $J_i$ as a real-valued function in the class $C( \R_+\times \R \times[0,T])$ such that for any $s\in [0,T]$, $\ell_i(\cdot,\cdot,s) $ is differentiable in $\R_+\times \R$, and the partial derivatives $\partial_y \ell_i$ and $\partial_a \ell_i$ are continuous with respect to $(y,a,s)$.
We also introduce a specific cost at the origin, that is, a function $\ell_*\in C([0,T])$.
We make the following assumptions.
\begin{description}
\item{[H1]}
There exists a positive constant $C_0$ and a real number $c_0$ such that $ \ell_i(y,a,s)\ge C_0|a|^p+c_0$, for every $i\in\{1,\dots,N\}$, $(y,a,s)\in \R_+\times \R\times[0,T]$ and $\ell_*(s)\geq c_0$ for any $s\in [0,T]$. There is no loss of generality in assuming $c_0=-C_0$, because if $c_0<-C_0$, we may always add the same positive constant to all the running and terminal costs. We therefore assume that there exists a positive constant $C_0$ such that 
for every  $(y,a,s)\in \R_+\times \R \times[0,T]$ and  $i\in\{1,\dots,N\}$,
 \begin{eqnarray}\label{H1}
 \ell_i(y,a,s)\ge C_0(|a|^p-1),\\
 \label{H1bis}
\ell_*(s)\geq -C_0.
\end{eqnarray}
\item{[H2]} For any $r>0$, there exists a modulus of continuity $\omega_r$ such that
\begin{equation*}
|\ell_i(y,a,t_1)-\ell_i(y,a,t_2)|\leq  \left(1+|a|^p\right) \omega_r (|t_1-t_2|),
\end{equation*}
for every $i\in\{1,\dots,N\}$, $y\in[0,r]$, $a\in\R$ and $t_1,t_2\in[0,T]$.
\item{[H3]}  For any $r>0$, there exists $\tilde c(r)\ge 0 $ such that 
\begin{equation*}%\label{H3}
\left|\partial_y\ell_i(y, a,s)\right|\leq \tilde c(r)(|a|^p+1)\quad\textrm{and }\quad   \left|\partial_a \ell_i(y, a,s)\right|\leq \tilde c(r)(|a|^{p-1}+1),
\end{equation*}
for every $i\in\{1,\dots,N\}$ and  $(y,a,s)\in [0,r]\times \R \times[0,T]$.
\item{[H4]}  $a\mapsto \ell_i(y,a,s)$ is strictly convex.
\end{description}

\noindent {\bf Notation}
When there is no ambiguity, for $J_i \ni x=\bar x e_i$ with $\bar x>0$ and for $\alpha=\bar \alpha e_i$,
we will sometimes write $\ell_i(x,\alpha,t)$ for $\ell_i(\bar x,\bar \alpha, t)$.

\begin{remark}\label{rmk:H1bis}
Assumption~$[H3]$ and the continuity of the running costs imply that for every positive number $r$, there exists a constant $\overline C_0(r)$ such that 
\begin{equation*}%\label{H3}
\left|\ell_i(y,a,s)\right|\leq \overline {C}_0 (r)(|a|^{p}+1)\qquad\forall(y,a,s)\in [0,r]\times \R \times[0,T] .
\end{equation*}
\end{remark}

\paragraph{Example}
Given $p>1$, $\kappa$ and $\lambda$ two bounded and continuous functions defined on $\R_+\times[0,T]$, ($\kappa$ being bounded from below by a positive constant), differentiable in $y$ with a bounded and continuous derivative,
$\ell_i: (y,a,s) \mapsto \kappa(y,s) |a|^p + \lambda(y,s)$ is an example of a function fulfilling Assumptions [H1] to [H4].

Before introducing the cost $J_t(x;y )$ associated with a trajectory $ y\in \Gamma_t[x]$, set $\bar y =|y|$, $\alpha= y'$, and $\bar\alpha=\bar y'$,  and define for any $\tau\in [0,T)$:
\begin{itemize}
\item   the  cost $\ell_\vertex$ for staying at the vertex $\vertex$,  obtained as the minimum of the costs
coming from the different edges and associated with zero velocity and of the  specific cost $\ell_*$: 
\begin{equation}
\label{eq:460}
\ell_\vertex(\tau)=\min\left\{\ell_*(\tau),\min_{i=1,\dots,N}\ell_i(0,0,\tau)\right\}
\end{equation}
\item the running cost
\begin{equation}\label{eq:462}
L(y(\tau), \alpha(\tau),\tau)=\sum_{i=1}^N \ell_i(\bar y(\tau),\bar \alpha(\tau), \tau) \car_{y(\tau)\in J_i\setminus\{\vertex\}}+\ell_\vertex(\tau)\car_{y(\tau)=\vertex}.
\end{equation}
\end{itemize}
Then, for $t\in [0,T]$, the cost associated with a trajectory $ y\in \Gamma_t[x]$ is 
\begin{equation}\label{eq:46}
J_t(x;y)=\int_t^T L\left(y(\tau), y'(\tau),\tau \right)\, d\tau+g(y(T)).
\end{equation}
It is convenient to drop the index when $t=0$, i.e. $J(x;y)=J_{0}(x;y)$.

The assumption on the structures of the costs at the origin in  \eqref{eq:461} and  \eqref{eq:460} will play an important role in the whole paper, in particular in the proof of the existence of optimal trajectories,  because they will imply the lower semi-continuity of  the cost functional in the space of admissible trajectories endowed with the weak topology  of $W^{1,p}(t,T;\R^N)$.

\paragraph{The value function}
The value of the optimal control problem is 
\begin{equation}\label{eq:4}
u(x,t)= \inf_{y\in \Gamma_t[x]}  J_t(x; y ).
\end{equation}
The set of optimal trajectories starting from $x$ at $t$ (i.e. which achieve $u(x,t)$) is 
\begin{equation}\label{eq:gamma_opt}
\Gamma^{\rm{opt}}_t[x]=\Bigl\{y\in\Gamma_t[x]\;:\; J_t(x;y)=\inf_{\hat y\in \Gamma_t[x]} J_t(x;\hat y)\Bigr\}.
\end{equation}
For simplicity, we drop the subscript when $t=0$: $\Gamma^{\rm{opt}}[x]=\Gamma^{\rm{opt}}_0[x]$.

\begin{remark}\label{rmk:bound_contr}
The value function $u$ is bounded locally uniformly in $\cG\times[0,T]$. Indeed, the trajectory associated with the control $\alpha\equiv 0$ is admissible and provides an upper bound on $u(x,t)$ depending possibly on $x$.  On the other hand, assumption [H1] and the boundedness of $g$ provide a lower bound on $u(x,t)$. 
\end{remark}
We skip the proof of the following lemma.
\begin{lemma}\label{lemma:u_g}
For every $x\in \cG$, $\lim_{t\to T^-} u(x,t)=g(x)$.
\end{lemma}

\subsection{Existence of optimal trajectories} \label{exist_opt_traj}
To our knowledge, Proposition \ref{prp:ex_OT} below on the existence of optimal trajectories is new with this degree of generality. While in \cite{AMMT3},  the existence of optimal trajectories was mostly a consequence of the lower semi-continuity of the cost w.r.t. the weak $L^2$-convergence of the controls, a similar semi-continuity property is not obvious in the present case. We need more involved arguments that we were not able to find in the available literature. Our strategy consists of seeing $\cG$ as a subset of~$\R^N$ (recall that $N$ is the number of edges; we already said in Section~\ref{sec:notations} that this does not bring any restriction to the hypotheses),  then extending in a suitable manner the cost to a functional acting on trajectories defined on the whole space~$\R^N$, and finally proving that the latter extension has properties that make it possible to invoke a general theorem from the monograph \cite{MR1020296}. Proposition \ref{prp:ex_OT} is a consequence of the following lemma, whose proof is postponed to Section \ref{proofs_of_technical_lemmas}.
\begin{lemma}\label{lem:LSC}
  Given $(x,t)\in \cG\times [0,T)$, the functional $j$ defined on $\Gamma_t[x]$ by
\begin{equation}\label{def:j}
  j(y)=\int_t^T  L(y(\tau),y'(\tau),\tau) d\tau,
\end{equation}
is %sequentially
$W^{1,p}$-weakly lower semi-continuous, i.e. if $(y^n)_n$ is a sequence in $\Gamma_t[x]$ such that
$y^n\rightharpoonup y$  weakly in $W^{1,p}(t,T;\R^N)$, (recall that  $y\in \Gamma_t[x]$, see Lemma \ref{compactness_lemma}), then  
\begin{equation}\label{LSC:1}
\liminf_n j(y^n)\geq j(y).
\end{equation}
\end{lemma}

\begin{proposition}\label{prp:ex_OT}
For each  $(x,t)\in \cG\times[0,T]$,  $\Gamma_t^{\rm{opt}}[x]\ne \emptyset$. In other words, there exists an optimal trajectory, i.e.  $y\in \Gamma_t[x]$ such that $u(x,t)=J_t(x; y)$.
Moreover, there exist positive numbers $M(|x|)$ and $R(|x|)$ such that if  $y\in \Gamma_t^{\rm{opt}}[x]$, then  $\|y\|_{W^{1,p} (t,T;\R^N)} \le M(|x|)$, $\|y\|_{C^{\frac {p-1}{p}}([t,T];\cG)}\leq M(|x|)$  and  $\max_{s\in [t,T]} |y(s)| \le R(|x|)$.
\end{proposition}
\begin{proof}
For $t=T$, there is nothing to do. Fix now $(x,t)\in\cG\times[0,T)$, consider a minimizing sequence $y^n\in \Gamma_t[x]$, i.e. such that $\lim_{n\to\infty} J_t(x; y^n)=u(x,t)$ and define $\alpha^n$ as the weak time derivative of $y^n$. Thanks to the assumptions on the costs, standard calculus yields $\|y^n\|_{W^{1,p}(t,T;\R^N)}\le C$ for a constant $C$ independent of $n$. Thanks to Lemma \ref{compactness_lemma}, there exists $y\in \Gamma_t[x]$ such that,  possibly after the extraction of a subsequence, $y_n$ converges to  $y$ in $W^{1,p} (t,T;\R^N)$  weakly and in $C([t,T];\cG) $.

To prove that $y$ is an optimal trajectory, i.e. that
$u(x,t)=J_t(x;y)$,  it is enough to establish that
\begin{equation}\label{claim2bis}
\liminf_nJ_t(x;y^n)\geq J_t(x;y).
\end{equation}
First, arguing exactly as in \cite[Prop 2.5]{AMMT3} and using in particular  \eqref{eq:461}, we obtain
%We first claim  that
\begin{equation}\label{part5}
\liminf_{n\to\infty}g(y^n(T))\geq  g(y(T)).
\end{equation}
Next, from  Lemma~\ref{lem:LSC}, we obtain that
\begin{equation}\label{LSC2}
\liminf_{n\to\infty} \int_t^T  L(y^n(\tau),\alpha^n(\tau),\tau) d\tau \ge  \int_t^T  L(y(\tau),y'(\tau),\tau) d\tau.
\end{equation}
Combining \eqref{part5} and \eqref{LSC2} yields \eqref{claim2bis} and the optimality of $y$.

We skip the easy proofs of the bounds on the optimal trajectories; they follow essentially from the assumptions on the costs.
\end{proof}

\begin{remark}\label{basic_facts_value}
\begin{enumerate}
\item It is easy to prove (by contradiction) that,
given $y\in \Gamma^{\rm{opt}}_t[x]$ and $\bar t \in [t,T]$, $y_{\mid[\bar t,T]}\in \Gamma^{\rm{opt}}_{\bar t}[y(\bar t)]$.
  \item 
From point 1, we deduce that for every $y\in \Gamma^{\rm{opt}}_t[x]$, there holds
\begin{equation*}
u(x,t)=u(y(\bar t),\bar t)+\int_t^{\bar t} L(y(\tau),y'(\tau),\tau)d\tau\qquad \forall \bar t\in[t,T].
\end{equation*}
\item 
Concatenating two optimal trajectories yields an optimal trajectory. More precisely, for every $y\in \Gamma^{\rm{opt}}_t[x]$, $\hat t\in(t,T)$ and $\hat y\in \Gamma^{\rm{opt}}_{\hat t}[y(\hat t)]$, the concatenation~$\tilde y$ of $y$ and $\hat y$ belongs to $\Gamma^{\rm{opt}}_t[x]$.
 \item 
The following statement can be seen as the ``converse'' of point 2.
If $u(x,t)=u(y_1(t_1),t_1)+\int_t^{t_1} L( y_1(\tau), y'_1(\tau),\tau) d\tau$ for some $t_1\in[t,T]$ and $y_1\in \Gamma_{t,t_1}[x]$, then there exists $y\in \Gamma^{\rm{opt}}_t[x]$ with $y_1=y$ on $(t,t_1)$.
\end{enumerate}
\end{remark}

\subsection{Proof of Lemma \ref{lem:LSC}}\label{proofs_of_technical_lemmas}
Recall that $N$ is the number of edges and that the $N$ unit vectors defining the orientations of the edges form the canonical basis of $\R^N$. This will play an important role in what follows. Recall also that the latter property is not a restrictive assumption, as explained in Section \ref{sec:notations}.

The strategy is as follows. In order to apply the powerful  \cite[Theorem 4.1.1]{MR1020296} from the calculus of variations, we need to extend the functional $j$ to non-admissible trajectories $y\in W^{1,p}([t,T];\R^N)$ with $y([t,T])\nsubseteq \cG$. The extended functional~${\mathcal J}$ is defined in~\eqref{def:cal_J} below, and involves the running cost~$\phi$ defined in~\eqref{eq:def_phi}. The lower semi-continuity of ${\mathcal J}$ with respect to the weak convergence in $W^{1,p}([t,T];\R^N)$ will be key to obtaining the existence of a minimizer. The function $(z,a,s)\mapsto \phi(z,a,s)$ takes its values in $(-\infty, +\infty]$, is convex in~$a$ and lower semi-continuous in~$(z,a)$ and such that ${\mathcal J}(y)$ coincides with~$j(y)$ if $y\in\Gamma_t[x]$ and is infinite otherwise  (note that no convexification is needed in the state variable $z$ due to the  relative compactness of the minimizing sequences of trajectories in $C([t,T];\cG)$). Those properties of $\phi$
stem in particular from the definition of $\phi(0,a,s)$, see  \eqref{eq:def_phi} and \eqref{def:Ell_0},  and from
the structural assumption in \eqref{eq:460}, which ensures the lower semi-continuity of the running cost at $(z,a)=(0,0)$.

Let the function $\phi:\R^N\times \R^N \times [0,T] \to (-\infty, +\infty]$ be defined by
\begin{equation}\label{eq:def_phi}
\phi(z,a,s)=\left\{\begin{array}{ll}
{\mathcal L}_0 (a,s)&\textrm{if }z=\vertex,\\
\ell(z,a,s)&\textrm{if }z\ne \vertex,
\end{array}\right.
\end{equation}
where
\begin{equation}\label{eq:def_l}
\ell(z,a,s)=\left\{ \begin{array}{ll} \ds
\ell_i(z\cdot e_i,a\cdot e_i,s)&\quad\textrm{if }\ds (z,a)\in J_i\times \R e_i \setminus\{(\vertex,0)\},\\
\ds \min_{i=1,\dots,N}\ell_i(0,0,s)&\quad\textrm{if }(z,a)=(\vertex,0),\\
+\infty &\quad\textrm{otherwise,}
\end{array}\right.
\end{equation}
and
\begin{equation} \label{def:Ell_0}
\mathcal L_0(a,s)=\conv\Bigl(\ell(\vertex,\cdot,s),\ell_\vertex(s)+\chi_{0}(\cdot)\Bigr) (a),
\end{equation}
in which the convex envelope is taken only with respect to the variable $a$, and $\chi_0$ is the characteristic function of $\{0\}$ in $\R^N$, i.e. it takes the value $0$ at $0$ and $+\infty$ elsewhere. \\
In the next lemma, $\mathbb B_{\R^N}$ and $\mathcal L_{[0,T]}$, respectively, stand for the Borel $\sigma$-algebra of $\R^N$
 and the Lebesgue $\sigma$-algebra of $[0,T]$.
\begin{lemma}\label{lemma:barL0}
The function~${\mathcal L}_0:\R^N\times[0,T]\to \R$  in~\eqref{def:Ell_0} has the following properties:
\begin{itemize}
\item[i)] $ {\mathcal L}_0( 0,s)=\ell_\vertex(s)$ for $s\in[0,T]$
\item[ii)] Given $\alpha\in\R^N$, the function ${\mathcal L}_0(\alpha,\cdot)$ is lower semi-continuous on $[0,T]$
\item[iii)]  ${\mathcal L}_0$ is $ \mathbb B_{\R^N} \otimes \mathcal L_{[0,T]}$-measurable.
\end{itemize}
\end{lemma}
\begin{proof}
The proof is rather technical and can be skipped on first reading. For this reason, it is postponed to the Appendix \ref{sec:app1}.
\end{proof}
\begin{lemma}\label{lemma:prp_cal_J}
 % The functional ${\mathcal J}$ and
  The function $\phi$ introduced in \eqref{eq:def_phi} has the following properties:
  %and respectively, have the following properties
\begin{itemize}
\item[i)] $\phi$ is $\mathbb B_{\R^N}\otimes \mathbb B_{\R^N} \otimes \mathcal L_{[0,T]}$-measurable
\item[ii)] Given $s\in[0,T]$, the function $\phi(\cdot,\cdot,s)$ is lower semi-continuous in  $\R^N\times\R^N$
\item[iii)] Given $(y,s)\in \R^N\times [0,T]$, the function $\phi(y,\cdot,s)$ is convex.
\end{itemize}
\end{lemma}
\begin{proof}
  Observe that 
\begin{equation*}
\phi(y,\alpha,s)=\left\{\begin{array}{ll}
+\infty & \textrm{if }y\notin\cG \textrm{ or } \textrm{if }y\in J_i\setminus\{\vertex\}\hbox{ and } \alpha\notin \R e_i\\
\ell_i(y\cdot e_i,\alpha\cdot e_i,s) &\textrm{if }y\in J_i\setminus\{\vertex\} \hbox{ and } \alpha\in \R e_i\\
{\mathcal L}_0(\alpha,s)&\textrm{if }y=\vertex.
\end{array}\right.
\end{equation*}
$i)$ Note first that the set $\{(y,\alpha,s)\in \R^N\times \R^N \times [0,T] \;:\;  \phi(y,\alpha, s)=+\infty\}$ is measurable.

Given $b \in \R$, we claim that the set $I=\{ (y,\alpha,s) \in \R^N\times \R^N \times [0,T] \;:\; \phi(y,\alpha,s)<b\}$ is
$\mathbb B_{\R^N}\otimes \mathbb B_{\R^N}\otimes \mathcal L_{[0,T]}$-measurable. Indeed, $I=\cup_{i=0}^N I_i$ where
$I_0=\{ (\vertex,\alpha,s) \;:\; \cL_0(\alpha,s)<b\}$ and $I_i=\{   (y,\alpha,s) \in J_i \setminus\{\vertex\} \times \R e_i \times [0,T] \;:\; \ell_i(y\cdot e_i,\alpha\cdot e_i,s)<b\}$ for $1\le i \le N$. The sets $I_i$ are clearly measurable if  $1\le i \le N$, and $I_0$ is measurable thanks to point $iii)$ in Lemma \ref{lemma:barL0}. Hence $I$ is measurable as the  union of $N+1$ measurable sets.
\\
$ii)$ Fix $s\in[0,T]$ and consider a  sequence $(y_n,\alpha_n)_n$ in $\R^N\times\R^N$ such that $(y_n,\alpha_n)\to (y,\alpha)$. We claim that $\liminf_n \phi(y_n,\alpha_n,s)\geq \phi(y,\alpha,s)$. To prove this, it suffices to focus on the case where $\liminf_n \phi(y_n,\alpha_n,s)\in \R$. Possibly after the extraction of a subsequence, we may assume that $\phi(y_n,\alpha_n,s)$ has a limit in $\R$ as $n\to \infty$. From the definition of $\phi$ and the closedness of $G$, this implies that $y\in \cG$. We make several cases:
\begin{itemize}
\item If  $y\in J_i\setminus\{\vertex\}$ for some  $i\in \{1,\dots, N\}$, then $y_n \in   J_i\setminus\{\vertex\}$ and $\alpha_n\in \R e_i$ for $n$ sufficiently large. Hence $\alpha\in \R e_i$ and the claim follows from the continuity of $\ell_i$.
\item If $y=\vertex$, then $\phi(y, \alpha, s)= \cL_0(\alpha, s)$. If $y_n=\vertex$ for $n$ sufficiently large, then the claim follows from the continuity of $\cL_0(\cdot,s)$.
\item If $y=\vertex$ and if $y_n\in \cG\setminus \{\vertex\}$ for $n$ large enough, then, after the extraction of a subsequence, we may assume that for some $i\in \{1,\dots, N\}$, $y_n\in J_i\setminus\{\vertex\}$ and $\alpha_n\in \R e_i$ for $n$ sufficiently large, so $\phi(y_n, \alpha_n, s)=\ell_i(y_n\cdot e_i, \alpha_n\cdot e_i, s)$. Since $\phi(\vertex, \alpha, s)=\cL_0(\alpha,s)\le 
  \ell_i(0, \alpha\cdot e_i, s)$, the  claim follows from the continuity of $\ell_i$.
\end{itemize}
$iii)$ The convexity of $\alpha\mapsto \phi(y,\alpha, s)$ follows from the convexity of $\ell_i$ for all $i\in \{1,\dots,N\}$ and of $\cL_0$ with respect to $\alpha$. 
\end{proof}

\begin{proof}[Proof of Lemma  \ref{lem:LSC}]
  Consider a new cost for trajectories $y\in  W^{1,p}([t,T];\R^N)$:
\begin{equation}\label{def:cal_J}
{\mathcal J}(y)=\int_t^T \phi\left(y(s),y'(s),s\right)\, ds,
\end{equation}
where $\phi$ is defined in \eqref{eq:def_phi}.

Observe that the three properties stated in Lemma~\ref{lemma:prp_cal_J} coincide  with the three assumptions in~\cite[Theorem 4.1.1]{MR1020296}. From Lemma~\ref{lemma:prp_cal_J} and \cite[Theorem 4.1.1]{MR1020296}, we find that the functional $\mathcal J$ is
% sequentially
lower semi-continuous with respect to the weak topology of $W^{1,1}([t,T];\R^N)$. This implies that $\mathcal J$ is lower semi-continuous with respect to the weak topology of $W^{1,p}([t,T];\R^N)$.

Using the same arguments as in the paragraph devoted to admissible trajectories in Section \ref{setting}, we see that for all $y\in\Gamma_t[x]$,  ${\mathcal J}(y)=j(y)$. Thus, $j$ is lower semi-continuous in $\Gamma_t[x]$ endowed with the weak topology of $W^{1,p}([t,T];\R^N)$, which is the desired conclusion.
\end{proof}

\subsection{First properties}

\begin{proposition}[Dynamic programming principle]\label{prp:prop_vf}
 For every $(x,t)\in\cG\times[0,T]$ and $\bar t\in[t,T]$, there holds 
\begin{equation*}\label{eq:DPP}
u(x,t)= \inf_{y\in \Gamma_{t,\bar t}[x]}\left\{
u(y(\bar t),\bar t)+\int_t^{\bar t} L(y(\tau),y'(\tau),\tau)d\tau
\right\}.
\end{equation*}
\end{proposition}
\begin{proof}
  We skip the proof for brevity.
\end{proof}
\begin{proposition}\label{prp:U_cont}  The function~$u$ is continuous in~$\cG\times[0,T)$.
\end{proposition}
\begin{proof}
We do not give details of the proof because it is quite similar to that of \cite[Proposition 2.14]{extended_version} (the reference \cite{extended_version} is a long version of \cite{AMMT3} that contains all the details).
Yet, we sketch the main arguments for completeness.
Consider $(x_1,t_1), (x_2,t_2)\in \cG\times[0,T_1)$, with $T_1<T$ and $\delta=d(x_1,x_2)<(T-T_1)/2$. There is no loss of generality in assuming that both~$x_1$ and~$x_2$ belong to the same edge, say~$e_1$. Consider $y_2\in\Gamma^{\rm{opt}}_{t_2}[x_2]$. The idea is to construct a trajectory $y_1\in \Gamma_{t_1}(x_1)$ that reaches $x_2$ at time $t_1+\delta$ with a constant velocity, then coincides with the composition of $y_2$ with an affine function of the time variable, more precisely
\begin{equation*}
y_1(s)=\left\{\begin{array}{ll}
x_1+\frac{x_2-x_1}{\d}(s-t_1) &\quad\textrm{for }s\in[t_1,t_1+\d],\\
y_2\left(\frac{T-t_2}{T-t_1-\d}s-\frac{T(\delta+t_1-t_2)}{T-t_1-\d}\right)
 &\quad\textrm{for }s\in[t_1+\d,T],
\end{array}\right.
\end{equation*}
so $y_1(T)=y_2(T)$.

Then, thanks to  assumption [H2] and Proposition \ref{prp:ex_OT}, easy computations yield the existence of  a modulus of continuity~$\tilde\omega$ such that 
\[u(x_1,t_1)-u(x_2,t_2)\leq J_{t_1}(x_1;y_1)-J_{t_2}(x_2;y_2)=\tilde\omega(|t_1-t_2|+\delta).
\]
Exchanging the roles of~$(x_1,t_1)$ and~$(x_2,t_2)$ then leads to the desired result.
\end{proof}

\begin{lemma}\label{lemma1}
Fix $x\in\cG$ and $y\in\Gamma[x]$. Consider a sequence $(x_n)_{n\in\N}$, $x_n\in\cG$ such that $ \d_n=  d(x_n,x)\to 0$ as $n\to \infty$.
There exists a sequence $(y_n)_{n\in\N}$, $y_n\in\Gamma[x_n]$,  such that, for every $n\in\N$,  $y_n(T)=y(T)$ and 
\begin{equation*}
\begin{array}[c]{rl}
(i)&\quad \ds \sup_{s\in [0,T]}d(y_n(s),y(s))\leq \d_n+\|y'\|_{L^p(0,T)}\d_n^{(p-1)/p} \\
(ii)&\quad\ds  \lim\limits_{n\to\infty}J(x_n;y_n) = J(x;y).
\end{array}
\label{eq:lemma1}
\end{equation*}
\end{lemma}

\begin{proof}%[Proof of Lemma~\ref{lemma1}]
Without any loss of generality, we may assume that (possibly  after the extraction a subsequence still denoted by $x_n$), the points $x$ and $x_n$ belong to the same edge ($J_1$ for instance).
Let us now introduce the admissible control 
\begin{equation*}
\alpha_n(s)=\left\{\begin{array}{ll}
e_1,\quad&\hbox{if } s\in[0,\delta_n] \hbox{ and }|x_n|\leq  |x|, \\
-e_1,\quad&\hbox{if } s\in[0,\delta_n] \hbox{ and } |x_n|> |x|, \\
\ds \frac{T}{T-\d_n}\alpha \left(T \frac{s-\d_n}{T-\d_n}\right),\quad &\hbox{if }s\in(\delta_n,T],
\end{array}\right.
\end{equation*}
and set $y_n(s)=x_n+\int_0^s \alpha_n(\tau)\, d\tau$. It is straightforward to check $(i)$ and $y_n(T)=y(T)$. Then 
\begin{multline*}\label{eq:rescal4}
J(x_n;y_n)-J(x;y)= \int_0^{\delta_n}L(y_n(s),y'_n(s),s)ds\\+\int_{\delta_n}^TL(y_n(s),y'_n(s),s)ds-\int_0^{T}L(y(s),y'(s),s)ds.
\end{multline*}
By the definition of $y_n$ and Remark~\ref{rmk:H1bis}, we see that $|\int_0^{\delta_n}L(y_n(s),y'_n(s),s)ds|\leq 2\delta_n \bar C_0(r)$ for $n$ sufficiently large, where $r=|x|+1$. Observe also that
$
\int_{\delta_n}^TL(y_n(s),y'_n(s),s)ds-\int_0^{T}L(y(s),y'(s),s)ds= A_n+B_n$, where 
\begin{eqnarray*}
A_n=\int_{0}^T\left[L\left(y(\tau), \frac{T}{T-\delta_n}y'(\tau),\frac{T-\delta_n}{T}\tau+\delta_n\right)-L(y(\tau),y'(\tau),\tau)\right]d\tau,\\
B_n=-\frac{\delta_n}{T}\int_{\delta_n}^TL\left(y(\tau), \frac{T}{T-\delta_n}y'(\tau),\frac{T-\delta_n}{T}\tau+\delta_n\right)d\tau.
\end{eqnarray*}
Both $A_n$ and $B_n$ vanish as $n\to\infty$. Collecting all the observations above yields $(ii)$.
\end{proof}
\begin{proposition}\label{prp:prop1}
Assume that the sequence $(x_n)_{n\in\N}$, $x_n\in\cG$, converges to $x\in \cG$ as $n\to\infty$. Consider $y_n\in\Gamma^{\rm{opt}}[x_n]$ for every $n\in\N$. Assume that, as $n\to\infty$, $y_n$  converges uniformly to $y$. Then $y\in \Gamma^{\rm{opt}}[x]$.
\end{proposition}
\begin{proof} The proof follows the same arguments of the proof of \cite[Proposition 2.19]{AMMT3} replacing \cite[Lemma 2.26]{AMMT3} with Lemma~\ref{lemma1}. We omit the details.
\end{proof}

\subsection{Necessary optimality  conditions}
\label{sec:euler_lagrange}
Below, we address situations in which it is possible to provide necessary optimality conditions for an optimal trajectory.
They  consist of a family of differential equations along with a condition at the horizon.
The following lemma deals with the Euler-Lagrange condition
in time intervals $(t_1,t_2)\subset (0,T)$ for which an optimal trajectory lies in the interior of a given edge.
\begin{lemma}\label{lemma:EL}
 Consider any $(x,t)\in \cG\times[0,T]$ and  $y\in \Gamma^{\rm{opt}}_t[x]$ such that, for some $t_1,t_2\in [t,T]$ and a given $i\in \{1,\dots, N\}$,  
$ y(s)\in J_i\setminus \{\vertex\}$ holds for all $s\in (t_1,t_2)$.
Then, setting $\bar y=|y|$, and $\bar \alpha=\bar y'$,
\begin{enumerate}
\item  $s\mapsto \partial_a \ell_i(\bar y(s),\bar \alpha(s),s)$ is absolutely continuous on $[t_1,t_2]$
\item $\frac {d}{ds}\partial_a \ell_i(\bar y(s),\bar \alpha(s),s)=\partial_y \ell_i(\bar y(s),\bar \alpha(s),s)$ 
  at a.e. $s\in[t_1,t_2]$
\item the function $y$ is Lipschitz continuous on $[t_1,t_2]$
 \item the  function $y$ is $C^1$  on $(t_1,t_2)$.
\end{enumerate}
\end{lemma}
\begin{remark}\label{rem:lip_constant}
A priori, the Lipschitz constant of $y$  appearing in Lemma \ref{lemma:EL}  may blow up as $t_2-t_1$ vanishes. We will later obtain a bound on the Lipschitz constant of $y$  that depends only on $T-t$ and $x$.
\end{remark}
\begin{proof}
Recall from Proposition \ref{prp:ex_OT} that there exists $R(|x|)>0$  such that $\bar y(s) \in [0,R(|x|)]$ for all $s \in [t_1,t_2]$.

Assume first that $ y(s)\in J_i\setminus \{\vertex\}$ for all $s\in [t_1,t_2]$.   Consider any $\bar \alpha_1\in L^p(t,T)$ such that
$\bar \alpha_1(s)=0$ for $s\notin (t_1,t_2)$  and $\int_{t_1}^{t_2}\bar \alpha_1 ds=0$. Set $\alpha_1(s)=\bar \alpha_1(s)e_i$.
For $h$ in an open interval containing $0$, the control $\alpha_h(\cdot)=\alpha(\cdot)+h\alpha_1(\cdot)$ is admissible for $(x,t)$ because $ \bar y|_{[t_1,t_2]}$ is bounded from below by a positive number.  Let $y_h$ denote the trajectory corresponding to the control $\alpha_h$. It is clear that $y_h(T)=y(T)$ and that $h\mapsto J_t(x; y_h)$ has a minimum at $h=0$.
 Therefore, for any $\bar \alpha_1$ as above and $\bar y_1(t)=\int_{t_1}^t \bar \alpha_1(s) ds$,
  \begin{displaymath}
    \int_0^T \partial_a \ell_i(\bar y(s),\bar \alpha(s),s)  \bar \alpha_1(s) +\partial_y \ell_i(\bar y(s),\bar \alpha(s),s)  \bar y_1(s) ds =0,
  \end{displaymath}
(note that the integral is meaningful thanks to assumption [H3] and the estimates on $y$). This yields that
\begin{equation}
\label{eq:opt_cond1}
\frac d {ds} \left(\partial_a \ell_i(\bar y(\cdot),\bar \alpha(\cdot),\cdot)\right)=\partial_y \ell_i(\bar y(\cdot),\bar \alpha(\cdot),\cdot),
\end{equation}
in the sense of distributions in $(t_1,t_2)$.

If $\bar y(t_1)=0$ or $\bar y(t_2)=0$, we can reproduce the argument in any subinterval $[\tilde t_1, \tilde t_2]$ strictly contained in $(t_1,t_2)$ and obtain again
 \eqref{eq:opt_cond1} in the sense of distributions in $(t_1,t_2)$.

Note that the right-hand side in  \eqref{eq:opt_cond1}  is a function in $L^1 (t_1,t_2)$ thanks to assumption [H3] and
that the function $s\mapsto \partial_a \ell_i(\bar y(s),\bar \alpha(s),s)$ is in $L^{\frac p {p-1}}(t_1,t_2)$.  From   \eqref{eq:opt_cond1},
  $s\mapsto \partial_a \ell_i(\bar y(s),\bar \alpha(s),s)$ belongs to $W^{1,1} (t_1,t_2)$, hence to $C([t_1,t_2])$ and $\|\partial_a \ell_i(\bar y(s),\bar \alpha(s),s)\|_{L^\infty(t_1,t_2)} \le C$, where, a priori, $C$ may blow up as $t_2-t_1$ vanishes.
  
This last point, the convexity of $\ell_i$ implied by assumption [H4], assumption [H1] and the bound on $\bar y$ yield 
\begin{equation*}
\begin{split}
C |\bar \alpha(s)|\ge \bar \alpha(s)  \partial_a \ell_i(\bar y(s),\bar \alpha(s),s) &\ge %\tilde C (|\bar \alpha (s) |^q -1)
\ell_i(\bar y(s),\bar \alpha(s),s) -\ell_i(\bar y(s),0,s) \\ &\ge  C_0(|\bar \alpha (s) |^p -1) -  \max_{|z|\le R(|x|),\; s\in [0,T]} \ell_i(|z|,0 ,s).    
\end{split}
\end{equation*}
Since $p>1$, this implies that for almost every $s\in (t_1,t_2)$, $|\bar \alpha(s)|\le A$ where $A$ is a constant that a priori blows up as  $t_2-t_1$ tends to $0$. Hence $y$ is Lipschitz continuous on $[t_1,t_2]$.

Finally, since $\partial_a \ell_i$  is strictly monotone with respect to its second argument (from  assumption [H4]) and  $\partial_a \ell_i$ is continuous in $\R_+\times \R\times [0,T]$,
and since $s\mapsto \bar y (s)$ is continuous, the continuity of $s\mapsto \frac {\partial\ell_i} {\partial a}(\bar y(s), \bar \alpha(s),s) $ in $(t_1,t_2)$ implies the continuity of $ \bar \alpha$ on $(t_1,t_2)$.  
\end{proof}
Next, we state a transversality condition for an optimal trajectory that stays in the interior of a given edge near the horizon $T$.
We need a further assumption on the terminal costs $g_i$:
\begin{description}
\item{[H5]}
  For all $i\in \{1,\dots, N\}$,   $g_i\in C^1(\R_+)$.
\end{description}

\begin{lemma}\label{lemma:trasver}
Under assumptions [H1] to [H5], consider any $(x,t)\in \cG\times[0,T)$ and any $y\in \Gamma^{\rm{opt}}_t[x]$ such that $y(T)\in J_i\setminus \{\vertex\}$. Then, with the same notation as in  Lemma~\ref{lemma:EL},
\begin{equation}\label{trasver_giugno}
\partial_a \ell_i(\bar y(T),\bar \alpha(T),T)=-\partial_x g_i(\bar y(T)).
\end{equation}
\end{lemma}
\begin{proof}
  The proof is similar to that of Lemma~\ref{lemma:EL}, see \cite[Theorem 6.2.4]{MR2041617}.
\end{proof}

\subsection{Lipschitz regularity of optimal trajectories}

We now aim at proving that for any $(x,t)\in \cG\times[0,T)$, any trajectory $y\in \Gamma^{\rm{opt}}_t[x]$ is Lipschitz continuous in $[t,T]$
with a Lipschitz constant that depends locally uniformly on $x$ and $T-t$ and that may explode as $t\to T$. 
The essential arguments are the Euler-Lagrange and the transversality conditions obtained in Section \ref{sec:euler_lagrange} and a key estimate on the initial velocity of an optimal trajectory, locally independent of the starting point, see Lemma \ref{lemma:casi_alphalim} below.\\
Recall from  Proposition~\ref{prp:ex_OT} that there exist $M(|x|)$ and $R(|x|)$, 
such that, for any $y\in \Gamma^{\rm{opt}}_t[x]$, 
\[\|y\|_{C^{\frac {p-1}{p}}([t,T];\cG)}\leq  M(|x|)  
%C_{\{|\xi|\leq |x|\}}
,\qquad \|y\|_{W^{1,p} (t,T;\R^N)} \le M(|x|)\quad\textrm{and }\max_{s\in [t,T]} |y(s)| \le R(|x|).
\]

\begin{theorem}\label{optcur_lip}
Under assumptions [H1] to [H5], for every $(x,t)\in \cG\times[0,T]$ and for every $y\in \Gamma^{\rm{opt}}_t[x]$, $y'$ belongs to $L^\infty(t,T;\R^N)$. Moreover, there exists a positive constant $V$ (depending only on $|x|$ and $T-t$, and that may explode as $t\to T$) such that
$\|y'\|_\infty \leq V$.
\end{theorem}
\begin{proof}
Consider $y\in \Gamma^{\rm{opt}}_t[x]$, set $\alpha=y'$, $\bar y =|y|\in W^{1,p}(t,T;\R_+)$ and $\bar\alpha= \bar y '\in L^p (t,T)$. 

Let us subdivide the interval $[t,T]$, to distinguish the times $s$ for which $y(s)\in J_i\setminus\{\vertex\}$, $i=1,\dots,N$, and $y(s)=\vertex$. More precisely, set
$I_0=\{s\in[t,T]:\; y(s)=\vertex\}$, and $I_i=\{s\in[t,T]:\; y(s)\in J_i\setminus\{\vertex\}\}$ for $i=1,\dots,N$.
Since $y$ is continuous, the set $I_0$ is closed and each $I_i$ can be written as the disjoint union of a (possibly infinite) family of
subintervals of $[t,T]$, relatively open in $[t,T]$.

We aim at bounding $\|\alpha\|_\infty$. For that, we consider the following different cases:

\begin{enumerate}
\item From Stampacchia theorem, $\alpha(s)=0$ for a.e. $s\in I_0$.
  
\item  Assume that, for some $i\in\{1,\dots,N\}$, we have either $(t_1, T]\subset I_i$ for some $t_1\in (t,T)$ or $[t,T]\subset I_i$ in which case we set $t_1=t$. This implies that $y(T)\in J_i\setminus\{\vertex\}$ and that $\alpha(s)=\bar \alpha(s)e_i$ and $y(s)=\bar y(s)e_i$ for $s\in [t_1,T]$. Hence \eqref{trasver_giugno} holds, together with \eqref{eq:opt_cond1} in $(t_1,T)$. From \eqref{eq:opt_cond1}, we deduce that $\| \frac d {ds} \partial_a \ell_i(\bar y(s),\bar \alpha(s),s)\|_{L^1(t_1,T)}$ is bounded by a constant which depends on $x$. Combined with \eqref{trasver_giugno}, this provides a bound on $\| \partial_a \ell_i(\bar y(s),\bar \alpha(s),s)\|_{L^{\infty}(t_1,T)}$ that depends only on $x$.  Then, arguing as in the end of the proof of Lemma \ref{lemma:EL}, we deduce that $\|\alpha\|_{L^\infty(t_1,T; \R^n)}$ is bounded by a constant that depends only on $x$.
\item Assume that, for some $t_1,t_2\in [t,T]$,  $y(t_1)=y(t_2)=\vertex$, and for some $i\in\{1,\dots,N\}$, $(t_1,t_2)\subset I_i$. Hence,  $\alpha(s)=\bar \alpha(s)e_i$ and $y(s)=\bar y(s)e_i$ for $s\in [t_1,t_2]$.
Since $y|_{(t_1,t_2)}$ is $C^1$, there exists $t^*\in (t_1,t_2)$ such that
$\bar \alpha (t^*)=0$. 
Combining this observation with \eqref{eq:opt_cond1}, we deduce that for all $s\in (t_1,t_2)$,
$|\partial_a \ell_i(\bar y(s),\bar \alpha(s),s)|\le   |\partial_a \ell_i(\bar y(t^*),0,t^*)| + \tilde c(R(|x|))\int_{t_1}^{t_2} (|\bar \alpha(s)|^p+1) ds $, ($\tilde c(\cdot)$ appears in Assumption [H3]). Hence,   $\| \partial_a \ell_i(\bar y(s),\bar \alpha(s),s)\|_{L^{\infty}(t_1,T)}$ is bounded by a constant that depends only on $x$.  Arguing as in the end of the proof of Lemma \ref{lemma:EL}, we deduce that $\|\alpha\|_{L^\infty(t_1,T; \R^n)}$ is bounded by a constant that depends only on $x$.

\item Assume that, for some $t_2\in (t,T]$,  $y(t_2)=\vertex$,   and for some $i\in\{1,\dots,N\}$, $[t,t_2)\subset I_i$. Hence,  $\alpha(s)=\bar \alpha(s)e_i$ and $y(s)=\bar y(s)e_i$ for $s\in [t,t_2]$.
We distinguish two subcases
\begin{enumerate}
\item  If there exists $s\in (t, t_2)$ such that $\bar y(s)>\bar y (t)$, then there exists $t^*\in (t, t_2)$ such that $\bar \alpha(t^* )=0$.  We may argue as in point $3$ and obtain the desired bound (depending only on $x$) on $\|\alpha\|_{L^\infty(t,t_2)}$.
\item  If $\bar \alpha(t)$ is negative, then we apply Lemma~\ref{lemma:casi_alphalim} below, which gives the bound on $\|\alpha\|_{L^\infty(t,t_2)}$.
\end{enumerate}
\end{enumerate}
\end{proof}
\begin{lemma}\label{lemma:casi_alphalim}
We make the same assumptions as in Theorem~\ref{optcur_lip}.
If, for some $t_1\in(t,T]$ and $i\in\{1,\dots,N\}$,  $y(t_1)=\vertex$,  $y(s)\in J_i\setminus\{\vertex\}$ for $s\in[t,t_1)$, and if $y'(t)\cdot e_i<0$, then
the function $[t,t_1) \ni s\mapsto y'(s) \cdot e_i$
is bounded from below by a constant that depends on $T-t$ and $|y(t)|$ and that blows up as $t\to T$.
\end{lemma}
 
\begin{proof}
We place ourselves in the same setting as in
Theorem~\ref{optcur_lip} and its proof, in particular for the definition of $x$ and $\alpha$.

The main part of the proof consists of obtaining a 
bound on the initial velocity, namely $|\alpha(t)\cdot e_i|$.
With the latter bound, arguments similar to the end of the proof of Lemma \ref{lemma:EL} will lead to the desired estimate.
The bound on the initial velocity will be found in an indirect manner: we will introduce a quantity that blows up as $|x|/ (T-t)$ tends to infinity, (see the right hand side of \eqref{eq:9999} below), and we will check that
\begin{itemize}
\item  if $|\alpha(t)\cdot e_i|$ is less than this quantity (namely \eqref{eq:9999} below), then we have obtained the desired bound
\item if \eqref{eq:9999} does not hold, then  $|\alpha(t)\cdot e_i|$ is bounded by a constant that does not depend on $\frac {|x|} {T-t}$ (this means that  \eqref{eq:9999} does not hold only if  $|x|=O(|T-t|)$ as  $t\to T$).
To prove this, we will construct suitable competitors of the optimal trajectory.
\end{itemize}
\noindent Several constants will be involved in this technical proof.
\begin{itemize}
    \item Some of them depend only on $\bar x =|x|$, namely:
 $\tilde V(\bar x)$, $\chi(\bar x)$, $k(\bar x)$, $m(\bar x)$, $c(\bar x) $, $V_*(\bar x) $, $r_*(\bar x)$ and $C^*(\bar x)$
 \item
 Other quantities will depend on $\bar  x$ and on $T-t$ and deteriorate as $t\to T$, see \eqref{eq:9999} below, but their dependence on $T-t $ will be  written explicitly.
\end{itemize}

Recall that $y'$ is continuous in  $[t,t_1)$. For any $s\in [t,t_1)$, set $y(s)=\bar y(s)e_i$ and $y'(s)=\bar \alpha(s)e_i$,
$\bar y(t)=\bar x$ and $\bar \alpha(t)=-\bar v$ with $\bar v>0$.
\\
We first claim that there exist constants $\tilde V(\bar x)$ and $\chi(\bar x)\ge 1$ such that, if $\bar v\ge \tilde V(\bar x)$, then
\begin{equation}\label{eq:1ago24}
\bar \alpha(s)\in[-\chi(\bar x) \bar v, -\bar v/\chi(\bar x)],\qquad \forall s\in [t,t_1).
\end{equation}
In fact, from the convexity of $\ell_i$,  
$a \partial_{a} \ell_i(\bar y, a, s) \ge \ell_i(\bar y, a, s) - \ell_i(\bar y, 0, s)$ for all $(\bar y, a,s) \in \R_+\times \R\times [0,T]$,
hence $a \partial_{a} \ell_i(\bar y, a, s)\ge C_0(|a|^p-1) - \ell_i(\bar y, 0, s)$ due to the assumption [H1]. Therefore, with $R(\bar x) $ defined in Proposition~\ref{prp:ex_OT}, 
\begin{displaymath}
  a \partial_{a} \ell_i(\bar y(s), a, s)\ge  C_0(|a|^p-1) -\max_{(z,s) \in [0,R(\bar x)]\times [0,T]}  \ell_i(z, 0, s) = C_0 |a|^p -  k(\bar x ),
\end{displaymath}
where $k(\bar x)=C_0 + \max_{(z,s) \in [0,R(\bar x)]\times [0,T]}  \ell_i(z, 0, s)$. Hence, also taking into account assumption~[H3], there exists a constant $m(\bar x)>0$  such that
for all $(s,a)\in [t,t_1)\times \R_+$, $\partial_{a} \ell_i(\bar y(s), a, s) \ge  C_0 |a|^{p-2}a - m(\bar x)$ and 
\begin{equation}\label{eq:20001}
\partial_{a} \ell_i(\bar y(s), -a, s) \le -C_0 |a|^{p-2}a +m(\bar x).
\end{equation}
On the other hand, let us set $K(\bar x)=\inf_{(z,s)\in [0,R(\bar x)]\times [0,T]}\partial_a \ell_i(z,0,s)$.
We also know from \eqref{eq:opt_cond1}, assumption~[H3] and Proposition~\ref{prp:ex_OT}, that for all $s\in [t,t_1)$,
\begin{equation}
  \label{eq:20002}
| \partial_{a} \ell_i(\bar y(s), \bar \alpha(s), s) - \partial_{a} \ell_i(\bar x, -\bar v, t) |\le c(\bar x).
\end{equation}
Combining   \eqref{eq:20002} with  \eqref{eq:20001}, we deduce that
$  \partial_{a} \ell_i(\bar y(s), \bar \alpha(s), s)\le -C_0 |\bar v|^{p-1} +m(\bar x)+c(\bar x)$.
 Therefore, from the strict convexity of $\ell_i$ in the variable $a$, if $\bar v$
 is sufficiently large so that $-C_0 |\bar v|^{p-1} +m(\bar x) +c(\bar x) < K(\bar x)$, then $ \bar \alpha(s)<0$ for all $s\in (t,t_1)$. Then, from assumption [H3], $-\tilde c(R(\bar x)) (|\bar \alpha(s)|^{p-1} +1) \le  -C_0 |\bar v|^{p-1} +m(\bar x) +c(\bar x)$  for all $s\in (t,t_1)$. Hence, if $\bar v$ is sufficiently large,
 then there holds $  \bar \alpha(s)\le -     \left (\frac {C_0}{2\tilde c(R(\bar x))}\right)^{\frac 1 {p-1}}  \bar v $ for all $s\in [t,t_1)$.
 Set $\chi(\bar x) = \max\left(  \left (\frac {2\tilde c(R(\bar x))} {C_0} \right)^{\frac 1 {p-1}},1\right)\ge 1 $. We have found that for $\bar v $ sufficiently large,  $  \bar \alpha(s)\le - \frac 1 {\chi(\bar x)} \bar v$ for all $s\in [t,t_1)$.
 
Similarly, $ \partial_{a} \ell_i(\bar y(s), \bar \alpha(s), s)\le -C_0  |\bar \alpha(s)|^{p-1} +m(\bar x)$. Then, from \eqref{eq:20002},
$  \partial_{a} \ell_i(\bar y(s), \bar \alpha(s), s)\ge  \partial_{a} \ell_i(\bar x, -\bar v, t)  -c(\bar x)\ge -\tilde c(R(\bar x)) (\bar v^{p-1} +1) -c(\bar x)$. 
 We deduce that for $\bar v$ large enough,  $  \bar \alpha(s)\ge -\chi(\bar x)\bar v $ for all $s\in [t,t_1)$. In conclusion, there exists a constant $\tilde V(\bar x)$ such that, if $\bar v\ge \tilde V(x)$, then \eqref{eq:1ago24} holds.

Observe that  \eqref{eq:1ago24}  implies that $\bar x-\chi(\bar x) \bar v(t_1-t)\leq 0\leq \bar x-\bar v(t_1-t)/\chi(\bar x)$, which also reads
\begin{equation}\label{eq:1ago24b}
t+   \chi^{-1}(\bar x)  \frac{\bar x}{ \bar v}\leq t_1\leq t+   \chi(\bar x) \frac{ \bar x}{\bar v}.
\end{equation}
Then
\begin{eqnarray*}
\int_t^{t_1}\ell_i(\bar y(s),\bar \alpha(s),s) ds&\geq& C_0\int_t^{t_1}\left(|\bar \alpha(s)|^p-1\right) ds
\geq C_0   \chi^{-1}(\bar x)  \frac{\bar x}{ \bar v}   \left(\frac{\bar v^p}{\chi^p(\bar x)}-1\right).
\end{eqnarray*}
Hence,  if 
$
\bar v \ge \max\left (\tilde V(\bar x), 2^{\frac 1 p}\chi(\bar x)\right)$, then
\begin{equation}\label{eq:1ago24c}
\int_t^{t_1}\ell_i(\bar y(s),\bar \alpha(s),s)\, ds \geq\frac {C_0}{2\chi^{p+1}(\bar x)}     \bar x \bar v^{p-1}.
\end{equation}

We are now going to find estimates on $\bar v$ by proposing
suitable competitors for the optimal control problem defining $u(x,t)$.
Let us first introduce some parameters that will arise in the arguments below:
\begin{itemize}
\item $V_*(\bar x)$ is the bound on $|y'|$ established in $(1)$-$(2)$-$(3)$-$(4.a)$ in the proof of Theorem~\ref{optcur_lip}.
\item $r_*(\bar x) =\max(\bar x, 2V_*(\bar x) T)$
\item $C^*(\bar x) = \max \left( \chi(\bar x) , \left (2 \chi^{p+1}(\bar x) \frac {\overline C_0(r_*(\bar x))}{C_0}\right)^{\frac 1 {p-1}} +1\right) $, where $C_0$ is the constant arising in assumption [H1] and $\overline{C}_0(\cdot)$ appears in Remark \ref{rmk:H1bis}.
  \end{itemize}
If 
\begin{equation}\label{eq:9999}
\bar v\le \max\left(  \tilde V(\bar x), 2^{\frac 1 p}\chi(\bar x),\chi(\bar x) \frac { \bar x }{T-t}, 2C^*(\bar x)\frac { \bar x }{T-t}\right),
\end{equation}
there is nothing to do. We are left with estimating $\bar v$ in the case where
\begin{equation}\label{eq:10000}
\bar v\ge \max\left( \tilde V(\bar x), 2^{\frac 1 p}\chi(\bar x), \chi(\bar x) \frac { \bar x }{T-t}, 2C^*(\bar x) \frac { \bar x }{T-t}\right).
\end{equation}
From \eqref{eq:1ago24b}, we deduce  that $t_1\le t+ C^*(\bar x) \frac {\bar x}{\bar v}\le T$   because   $C^*(\bar x)\ge \chi(\bar x)$  and  $\bar v \ge  2C^*(\bar x) \frac { \bar x }{T-t}$.
We can therefore construct $y_1\in \Gamma_{t, t+ C^*(\bar x) \frac {\bar x}{\bar v}}[x]$ by setting
$y_1'(s)=-\frac{\bar v}{C^*(\bar x)}e_i$ for $s\in \left[t,t+ C^*(\bar x) \frac {\bar x}{\bar v}\right]$. We see that $y_1\left(t+ C^*(\bar x) \frac {\bar x}{\bar v}\right)=\vertex$
and that $y$ reaches $\vertex$ before $y_1$.

The arguments in the following differ according to the behaviour of $y$ after $t_1$.

\noindent{\it Case A: $y\left(t+ C^*(\bar x) \frac {\bar x}{\bar v}\right)=\vertex$.}
Consider the admissible trajectory $y_A$ that coincides with $y_1$ in $ [t,t+  C^*(\bar x) \bar x/\bar v]$ and with $y$ in $[t+\bar x C^*(\bar x)/\bar v,T]$.
Set $\alpha_A=y_A'$ and define $\bar y_A(s)=|y_A(s)|$ and $\bar\alpha_A=\bar y_A'$. The optimality of $y$ entails
\begin{multline}\label{caseA}
0\leq J_t(x;y_A)-J_t(x;y)=
\int_{t}^{t+  C^*(\bar x) \bar x/\bar v}\ell_i(\bar y_A(s),\bar \alpha_A(s),s) ds\\-\int_t^{t_1}\ell_i(\bar y(s),\bar\alpha(s),s)ds-\int_{t_1}^{t+
 C^*(\bar x) \bar x/\bar v}L(y(s),\alpha(s),s)ds.
\end{multline}
From Remark \ref{rmk:H1bis}, we deduce that
\begin{equation*}\label{eq:26ago}
\int_{t}^{t+  C^*(\bar x) \bar x/\bar v}\ell_i(\bar y_A(s),\bar \alpha_A(s),s) ds\leq \overline C_0(r_*(\bar x))\bar x
\left(\frac{\bar v^{p-1}}{(C^*(\bar x))^{p-1}}+\frac{C^*(\bar x)}{\bar v}\right),
\end{equation*}
while assumption [H1] yields
\begin{equation*}
-\int_{t_1}^{t+  C^*(\bar x) \bar x/\bar v}L(y(s),\alpha(s),s)ds\leq C_0 C^*(\bar x) \frac{\bar x }{\bar v}.
\end{equation*}
From \eqref{caseA}, \eqref{eq:1ago24c} and the latter two inequalities, we obtain
\begin{displaymath}
0\leq \overline C_0(r_*(\bar x))\bar x\left(\frac{\bar v^{p-1}}{(C^*(\bar x))^{p-1}}+\frac{C^*(\bar x)}{\bar v}\right)+C_0  C^*(\bar x) \frac{\bar x}{\bar v}-
\frac {C_0}{2\chi^{p
+1}(\bar x)}     \bar x \bar v^{p-1} 
\end{displaymath}
i.e.
\begin{displaymath}
\bar v \le   M_A:=     \left(\frac { C^*(\bar x) (\overline C_0(r_*(\bar x)) +C_0 )}
{ \frac{C_0}{2\chi^{p+1}(\bar x)}\ - \frac{\overline C_0(r_*(\bar x))}{(C^*(\bar x))^{p-1}}}\right)^{\frac 1 p}.
\end{displaymath}
Note that the denominator of $M_A$ is positive from the definition of $C^*(\bar x)$ and that $M_A$ depends on $\bar x$ and not on $T-t$.

\noindent{\it Case B: $y(t+  C^*(\bar x) \bar x/\bar v)\in J_j\setminus\{\vertex\}$ for some $j\in \{1,\dots,N\}$.}
For brevity, we do not distinguish the case where $j=i$, although a sharper estimate on $\bar v$ can be found easily in this case.
From points $(1)$, $(2)$ and $(3)$ in the proof of Theorem~\ref{optcur_lip}, we know that $\|y'\|_{L^\infty(t_1,T)}\leq V_*(\bar x)$.
Let the  auxiliary trajectory  $y_2$  coincide with $y_1$ in $[t,t+C^*(\bar x) \bar x /\bar v]$ and satisfy  $y'_2(s)=2V_*(\bar x) e_j $ in  $ [t+  C^*(\bar x)\bar x/\bar v,
t+2C^*(\bar x) \bar x /\bar v]$. Note that $t+2C^*(\bar x) \bar x /\bar v\le T$  from the definition of $C^*(\bar x)$ and that the trajectory $y_2$ belongs to $\Gamma_{t, t+ 2C^*(\bar x) \frac {\bar x}{\bar v}}[x]$. Set $\alpha_2=y_2'$, and define $\bar y_2(s)= |y_2(s)|$ and $\bar \alpha_2= \bar y_2'$.
Observe next that, at least in a right neighbourhood of $t+ C^*(\bar x)\bar x /\bar v$, both $y(s)$ and $y_2(s)$ belong to $J_j\setminus\{\vertex\}$ and $\bar \alpha(s)-\bar \alpha_2(s)\leq -V_*(\bar x)$. 

We claim that there exists $t_*\in [t+C^*(\bar x) \bar x /\bar v, t+2C^*(\bar x) \bar x /\bar v]$ such that $\bar y(t_*)=\bar y_2(t_*)$. Indeed, if $y$ does not stay in $J_j$ in the interval $[t+C^*(\bar x) \bar x /\bar v, t+2C^*(\bar x) \bar x /\bar v]$, the two trajectories $y$ and $y_2$ must cross each other before $ t+2C^*(\bar x) \bar x /\bar v$ and the claim is proved.  If $y(s)\in J_j$ for all
$[t+C^*(\bar x) \bar x /\bar v, t+2C^*(\bar x) \bar x /\bar v]$, we argue as follows: 
since $y(t_1)=\vertex$, $
\bar y(t+ C^*(\bar x) \bar x /\bar v)\leq V_*(\bar x) \;(t+C^*(\bar x) \bar x /\bar v-t_1)\leq V_*(\bar x) C_* \bar x /\bar v$ and  $\bar y (s) -\bar y_2(s) \le V_* C_* \bar x /\bar v -V_* (s-t -C^*(\bar x) \bar x /\bar v)$. The right-hand side in the latter inequality is negative at $s=  t+2C^*(\bar x) \bar x /\bar v$,  which proves the claim. Moreover, from the definition of $r_*(\bar x)$,  it is easy to check that  $|y_2(s)| \le r_*(\bar x)$ for all $s\in [t, t_*]$.

We can now construct the desired competitor $y_B\in \Gamma_t[x]$ as follows: $y_B(s)= y_2(s)$ for all $s\in [t, t_*]$ and $y_B( s)= y(s)$ in $[t_*,T]$. From the optimality of $y$, we deduce that, with the now usual notations, 
\begin{eqnarray}\notag
0&\leq &J_t(x;y_B)-J_t(x;y)\\ \notag&=&
\int_{t}^{t+ \bar xC_*/\bar v}\ell_i(\bar y_B(s),\bar \alpha_B(s),s) ds+\int_{t+ \bar xC_*/\bar v}^{t_*}\ell_j(\bar y_B(s),\bar \alpha_B(s),s) ds\\\label{caseB}&&-\int_t^{t_1}\ell_i(\bar y(s),\bar\alpha(s),s)ds-\int_{t_1}^{t_*}L(y(s),\alpha(s),s)ds,
\end{eqnarray}
We already know that $ \int_{t}^{t+  C^*(\bar x) \bar x/\bar v}\ell_i(\bar y_B(s),\bar \alpha_B(s),s) ds\leq \overline C_0(r_*(\bar x))\bar x
\left(\frac{\bar v^{p-1}}{(C^*(\bar x))^{p-1}}+\frac{C^*(\bar x)}{\bar v}\right)$. %, see \eqref{eq:26ago}. 
Remark \ref{rmk:H1bis} implies $ \int_{t+  C^*(\bar x) \bar x/\bar v}^{t_*} \ell_j(\bar y_B(s),\bar \alpha_B(s),s) ds \le \overline C_0(r_*(\bar x))  ( 2^p V_*^p +1) C^*(\bar x) \bar x/\bar v $, because $t_*\le t+2C^*(\bar x) \bar x/\bar v$. Finally,
$-\int_{t_1}^{t_*}L(y(s),\alpha(s),s)ds\le C_0 (t_*-t_1)\le 2C_0C^*(\bar x)  \bar x/\bar v $. Combining \eqref{eq:1ago24c} and the last three estimates, \eqref{caseB} yields:
\begin{displaymath}
\bar v \le   M_B:=  \left(\frac { C^*(\bar x) (\overline C_0(r_*(\bar x))(2^pV_*^p+2) +2C_0 )}
{ \frac{C_0}{2\chi^{p+1}(\bar x)}\ - \frac{\overline C_0(r_*(\bar x))}{(C^*(\bar x))^{p-1}}}\right)^{\frac 1 p}.
\end{displaymath}
As in Case A, the denominator of $M_B$ is positive and $M_B$ depends on $\bar x$ and not on $T-t$.

We have bounded $|\bar v|$  by a constant that depends on $T-t$ and $\bar x$, and we conclude arguing as in the end of the proof of Lemma \ref{lemma:EL}.
\end{proof}

\subsection{Properties of the value function}
\label{sec:prop_u}

\subsubsection{Lipschitz continuity}
\label{sec:u_lip}

\begin{proposition}\label{prop:loclip}
Under assumptions [H1] to [H5], the value function is locally Lipschitz continuous in $\cG\times [0,T)$. 
\end{proposition}
\begin{remark}
Note that we do not suppose that the costs $\ell_i$ are Lipschitz continuous with respect to time.
\end{remark}
\begin{proof}
For brevity, the proof is only sketched because it is like \cite[Proof of Proposition~2.21]{AMMT3}  and borrows some ideas from \cite[Proposition 4.1]{CCC1}.

Set $\cG_{\delta}=\{ x\in \cG:  d(x,0)\le \delta\}$. We first prove that $u(\cdot,t)$ is locally Lipschitz continuous in $J_i\setminus\{0\}$ locally uniformly with respect to $t\in[0,T)$. More precisely, having fixed $T_1\in(0,T)$ and $R>0$, we wish to prove that  for any $t\in[0,T_1]$, $x_0= \bar x_0 e_i \in J_i \cap \cG_R \setminus\{0\}$ and $r>0$ sufficiently small, the function~$u(\cdot,t)$ is Lipschitz continuous on $(\bar x_0-r,\bar x_0+r)e_i$ with a Lipschitz constant that depends only on the parameters of the problem and on $T_1$ and $R$ (it is independent of $\bar x_0$, $r$, $t$ and $i$).

For that, fix  some $r$, $0<r< \bar x_0/4$ and observe that $(\bar x_0-4r,\bar x_0+4r)e_i\subset J_i\setminus\{0\}$. Consider $\bar x,\bar x_1\in (\bar x_0-r,\bar x_0+r)$, with $\bar x\ne\bar x_1$ and $|\bar x-\bar x_1|\leq 2(T-t)V$, where $V$ is the constant found in Theorem~\ref{optcur_lip} for the set $\cG_{5R/4} \times[0,T_1]$. Set $x=\bar x e_i$, $x_1=\bar x_1 e_i$ and $\tau=\frac{|\bar x-\bar x_1|}{ 2 V}$. For $y\in \Gamma^{\rm{opt}}_t[x]$,  let $y_1$ be the trajectory starting at $x_1$ at time~$t$ such that 
\begin{equation*}
y'_1(s)=y'(s)+(x-x_1)/\tau \quad \textrm{for }s\in[t,t+\tau],\qquad
y'_1(s)=y'(s) \quad \textrm{for }s\in(t+\tau,T].
\end{equation*}
From Theorem~\ref{optcur_lip}, $y$ does not reach the origin~$0$ before time $t+\frac{3r}{ V}$.
On the other hand, $\tau\leq \frac{r}{V}$.
%$\tau\leq \frac{r}{2V}$. 
Hence, in the time interval $(t,t+\tau)$, $y$ stays in $J_i\setminus\{0\}$. 

It is clear that  $y(\cdot)= y_1(\cdot)$ in $(t+\tau,T]$ and that, in the interval $(t,t+\tau)$, both $y$ and $y_1$ stay in $J_i\setminus\{0\}$. The same arguments as in ~\cite[Proof of Proposition~2.21]{AMMT3} yield that
$y_1\in \Gamma_t[x_1]$,  $d(y(\cdot),y_1(\cdot))\leq |\bar x-\bar x_1|$ in $(t,T]$.

By definition of $u$, following standard calculations and taking advantage of assumption~[H3] and Proposition~\ref{prp:ex_OT}, we obtain
\begin{eqnarray*}
u(x_1,t)-u(x,t)&\leq& 
\int_t^{t+\tau}\left|L(y_1(s),y'_1(s),s) - L(y_1(s),y'(s),s)\right|\, ds\\
&&\quad+\int_t^{t+\tau}\left|L(y_1(s),y'(s),s) - L(y(s),y'(s),s)\right|\, ds\\
&\leq & K |\bar x-\bar x_1|
\end{eqnarray*}
where $K$ is a constant depending only on~$|x_0|$.
Exchanging the roles of $x$ and $x_1$ yields the desired Lipschitz continuity locally uniform in~$x$.

The remainder of the proof contains the same arguments as in \cite[Proof of Proposition~2.21]{AMMT3}, and we skip it.
\end{proof}
\begin{remark}
If, furthermore, $g\in C(\cG)$, then the value function is locally Lipschitz continuous in $\cG\times [0,T]$. \end{remark}

%------------------------------------------------------------------------------------
\subsubsection{ The PDE satisfied by $u$ on  $\cG\times (0,T)$} 
\label{SHJB}
Below, we see that the value function is a viscosity solution
in a suitable sense
of an Hamilton-Jacobi equation in the network with a particular condition at the vertex, coming from the fact that, when standing at the vertex, the controller must take into account the discontinuities in the running costs occurring there. 
Recently, there has  been a growing research effort on Hamilton-Jacobi equations with discontinuous Hamiltonians and optimal control problems with discontinuities in  cost or dynamics. In particular, important progress has been made in the analysis of Hamilton-Jacobi equations posed on heterogeneous structures such as networks \cite{IM,MR3358634, MR3556345} (note that  \cite{MR3358634, MR3556345} are devoted to infinite horizon problems), {\sl booklet}-like geometries or {\sl multidimensional junctions} \cite{%barles2013bellman,
MR3690310, oudet2015equations}.
The first part of the book of G. Barles and E. Chasseigne, \cite{MR4704058} provides a state of the art on  optimal control problems with  discontinuities located on submanifolds of codimension one, and contains a chapter devoted to junctions.
In what follows, only assumptions [H1] to [H4] are needed.\\
\paragraph{Hamiltonians} For $x\in \R_+$, $i=1,\dots,N$, $t\in [0,T]$, $\bar p\in\R$, $p=(p_1,\dots,p_N) \in\R^N$, set
\begin{eqnarray*}
H_i(x,\bar p,t)&=&\sup_{ \zeta\in \R   
}\{-\bar p\zeta -\ell_i(x,\zeta,t)\},\\
H_i^{\downarrow}(x,\bar p,t)&=&\sup_{\zeta\ge 0
}\{-\bar p\zeta-\ell_i(x,\zeta,t)\}\\
H_\vertex(p,t)%&=&\max\Bigl(-\ell_*(t),\:\max_{i=1,\dots,N}\left\{-\ell_i(\vertex,0,t)\right\},\:\max_{i=1,\dots,N}\left\{H_i^{\downarrow}(\vertex,p_i,0)\right\}\Bigr)\\
&=&\max\Bigl(-\ell_\vertex(t),\:\max_{i=1,\dots,N}\left\{H_i^{\downarrow}(0,p_i,t)\right\}\Bigr).
\end{eqnarray*}
\paragraph{Hamilton-Jacobi equations on $\cG\times (0,T)$}
We are interested in the system of first-order PDEs on $\cG\times (0,T)$:
\begin{equation}\label{HJ}
\left\{\begin{array}{rcll}
-\partial_t u(x,t)+ H_i(|x|,D_xu,t)&=&0,&\qquad\textrm{if } x\in J_i\setminus\{\vertex\},\, 0< t<T,\\
-\partial_t u(\vertex,t)+ H_\vertex(D_xu,t)&=&0,&\qquad\textrm{if }  0< t<T,
\end{array}\right.
\end{equation}
where $D_xu(x,t) $ is defined in~\eqref{eq:def_deriv}.
 The contribution of $H_i^\downarrow$ to the Hamiltonian at the vertex (namely $H_0$) accounts for the trajectories that leave $0$ and enter the edge $J_i$, whereas $\ell_0(t)$ corresponds to the instantaneous cost of staying at the vertex. The quantity $-\ell_0(t)$ is seen as a flux limiter by Imbert and Monneau in
\cite{IM}. The following definition can be found in e.g. \cite{IM}:
\begin{definition} 
A function $u\in C(\cG\times(0,T))$ is a viscosity subsolution (resp. supersolution) of \eqref{HJ} if for every function $\varphi\in C^1(\cG\times(0,T))$ touching~$u$ from above (resp. below) at $(x,t)\in\cG\times(0, T)$, there holds
\begin{equation*} 
\begin{array}{rcll}
-\partial_t \varphi(x,t)+ H_i(|x|,D_x\varphi,t) &\leq& 0\quad (\textrm{resp. }\geq 0)&\qquad\textrm{if } x\in J_i\setminus\{\vertex\},\\
-\partial_t \varphi(\vertex,t)+ H_\vertex(D_x\varphi,t)&\leq& 0\quad (\textrm{resp. }\geq 0)&\qquad\textrm{if } x=\vertex.
\end{array}
\end{equation*}
A function $u\in C(\cG\times(0,T))$ is a viscosity solution of \eqref{HJ} if it is a viscosity subsolution and a viscosity supersolution of \eqref{HJ}.
\end{definition}
The proof of the following theorem can be found in \cite{extended_version}:
\begin{theorem}\label{thm:HJ}
The value function $u$ defined in~\eqref{eq:4} is a viscosity solution of \eqref{HJ} in $\cG\times (0,T)$. Moreover, for all $x\in \cG$, $t\mapsto u(x,t)$ is continuous in $ [0,T]$ and $u(x,T)=g(x)$.
\end{theorem}
If $g\in C_b(\cG)$, then there is a unique viscosity solution~$u$ to~\eqref{HJ} such that $u\in C_b(\cG\times[0,T])$ and $u(\cdot,T)=g$; see \cite{IM} for a proof. We have not found any literature on
uniqueness when g is not continuous and plan to address this topic in the future.
\subsubsection{Local semi-concavity of the value function away from the vertex}\label{sec:semi_concavity}

The following lemma will be useful in Subsection \ref{sec:regularity_u}, for studying the regularity of $u$ in the interior of the edges. 
\begin{lemma}\label{lemma:CC2_3.1}
Consider $x\in J_i\setminus\{\vertex\}$ for some $i=1,\dots,N$, $t\in[0,T)$ and $y\in\Gamma^{\rm{opt}}_t[x]$. Set $x=\bar x e_i$.  Under assumptions [H1] to [H5], for any $h\in \R$ such that $|h|<\bar x /2$,
\begin{equation}\label{eq:super_dif}
u(x+he_i,t)-u(x,t)+  h\partial_a \ell_i (\bar x ,  y'(t)\cdot e_i,t) \leq |h|\omega(|h|)
\end{equation}
where $\omega(\cdot)$ is a modulus of continuity which  is uniform w.r.t. $x$ in any compact contained in $J_i\setminus\{0\}$ and deteriorates as $t\to T$.
In particular,
$-\partial_a \ell_i (\bar x ,  y'(t)\cdot e_i,t)$ belongs to $D^+_x u(x,t)$, the superdifferential of $u(\cdot,t)$ at $x$, defined by
\begin{equation} \label{eq:super_diff_u_x}
 %-\partial_a \ell_i (\bar x ,  \alpha(t),t)\in
\ds D^+_{x}u(x,t)= \{q\in \R: \limsup_{J_i\ni y\to x}\frac{u(y,t)-u(x,t)-q  ( y- x)\cdot e_i}{| y- x|}\leq 0\}.\end{equation}
Finally, for any $t\in [0,T)$, $u(\cdot, t)$ is locally semi-concave with respect to $x\in \cG\setminus\{0\}$
in the sense of \cite[Definition 2.1.1]{MR2041617}.
\end{lemma}
\begin{proof}
Set $\bar y =|y|$. Recall from 
Proposition \ref{prp:ex_OT} and Theorem~\ref{optcur_lip} that $\|y\|_\infty\leq M $ and $\|y'\|_\infty\leq V$, where $M$ is a constant that depends only on $|x|$ while $V$ is a positive constant that depends on $|x|$ and on $T-t$ and deteriorates as $t\to T$. Take $h\in \R $ with $|h|< \bar x/2 $, and set $t_*=\left(\frac{T-t}{2}\right)\wedge \left(\frac{\bar x }{2V}\right)$.
Observe that $y(s)=\bar y(s) e_i$ for  $s\in [t,t+t_*]$.
The trajectory $y_h$ defined by
$y_h(s)= \bar y_h(s) e_i$ with $  \bar y_h(s)= \bar y(s)+ h ( t+t_* -s)/ t_* $ for $s\in[t,t+t_*]$, and $y_h(s)=y(s)$ for $s\in[t+t_*,T]$, is an admissible competitor for $u(x+he_i,t)$. Furthermore, 
it is easy to check that for all $s\in [t,t+t_*]$, $(\bar y_h(s),  \bar y_h'(s))\in [0,2\bar x] \times [-V-\frac {\bar x}{2t_*}, V+ \frac {\bar x}{2t_*}]$. Then,
\begin{equation*} u(x+he_i, t) -u (x,t)\le I +II,
\end{equation*}
where
\begin{equation*}
\begin{split}
I= &  -\frac {h}{t_*}\int_t^{t+t_*} \int_0^1 \partial_a \ell_i\left(\bar y_h(s), \bar y'(s) -\frac {\theta h}{t_*},s\right) d\theta ds,\\
II=& \frac {h}{t_*} \int_t^{t+t_*} \int_0^1   (t+t_* -s) \partial_y \ell_i \left( \bar y(s) +\theta (y_h(s)-y(s)), \bar y'(s) ,s\right) d\theta ds.
    \end{split}
  \end{equation*}
From the definition of $t_*$, assumptions [H1] and [H3] and the continuity of $\partial _a \ell_i$, there exists a modulus of continuity $\omega_1$ (uniform w.r.t. $x$ in any compact of $\cG$ contained in $J_i\setminus\{0\}$ and that deteriorates as $t\to T$) such that
\begin{equation*}
\ds I \le - \frac {h}{t_*}\int_t^{t+t_*} \partial_a \ell_i(\bar y(s), \bar y'(s) ,s)  ds +|h|\omega_1(|h|) .
\end{equation*}
From Lemma \ref{lemma:EL}, we can integrate by part in the latter formula and find 
\begin{equation*}
I\le-h \partial_a \ell_i(\bar y(t), \bar y'(t) ,t)-
\frac {h}{t_*} \int_t^{t+t_*}   (t+t_* -s) \partial_y \ell_i(\bar y(s), \bar y'(s) ,s)   ds +|h|\omega_1(|h|).
\end{equation*}
 Similarly, from the definition of $t_*$, assumption [H3] and  the continuity of $\partial_y \ell_i$, there exists a modulus of continuity $\omega_2$ (uniform w.r.t. $x$ in any compact of $\cG$ contained in $J_i\setminus\{0\}$ and that deteriorates as $t\to T$)
 such that
\begin{displaymath}
\begin{array}[c]{ll} &\ds II-
\frac {h}{t_*} \int_t^{t+t_*}  (t+t_* -s) \partial_y \ell_i(\bar y(s), \bar y'(s) ,s)   ds
\\ =&\ds \frac {h}{t_*} \int_t^{t+t_*}  (t+t_* -s) \int_0^1 \left(\partial_y \ell_i(\bar y(s) +\theta (\bar y_h(s) -\bar y(s)), \bar y'(s) ,s) - 
\partial_y  \ell_i(\bar y(s), \bar y'(s) ,s) \right) d\theta ds
\\\le & |h|\omega_2(|h|).
\end{array}
\end{displaymath}
Combining all the information above, we obtain
\eqref{eq:super_dif}.

Finally, from \cite[Theorem 2.1.10]{MR2041617},
the desired local semi-concavity w.r.t. $x\in \cG\setminus\{0\}$ is obtained  by
taking $h$ and $-h$ in \eqref{eq:super_dif} and summing.
\end{proof}
The following proposition stems from \cite[Proposition 5.3.8]{MR2041617}, see also \cite{MR986450}.
\begin{proposition}\label{prop:semiconc_loc}
Under assumptions [H1] to [H5] and if the functions $t\mapsto \ell_i(x,a,t)$ are Lipschitz-continuous with respect to $t$ with a Lipschitz constant that is locally uniform in $(x,a)$, the function $u$ is locally semi-concave in $ (J_i\setminus\{\vertex\}) \times (0,T)$ in the sense of \cite[Definition 2.1.1]{MR2041617}.
\end{proposition}

\begin{remark}\label{SCS_with_linear_modulus}
Under the additional assumption $[H+]$: 
\begin{description}
\item{[H+]}
$\;$ For all $s\in [0,T]$,  $ \ell_i (\cdot,\cdot,s) \in W^{2,\infty}_{\rm loc} ([0,\infty)\times \R)$, and for any compact subset $K$ of $[0,\infty)\times \R$, there exists a constant $M_K$ such that $\|\ell_i(\cdot,\cdot, s)\|_{W^{2,\infty}(K)} \le M_K$ for all $s\in [0,T]$,
\end{description}
it is possible to obtain more precise information on $u(\cdot,t)$. 

\medskip

Consider $x\in J_i\setminus\{\vertex\}$ for some $i=1,\dots,N$, $t\in[0,T)$ and $y\in\Gamma^{\rm{opt}}_t[x]$. Set $x=\bar x e_i$.  Under assumptions [H1] to [H5] and [H+], there exists a constant $C$ (depending on $|\bar x|$ and $T-t$) such that for every $h\in \R$ sufficiently small,
\begin{equation*}
u(x+he_i,t)-u(x,t)+  h\partial_a \ell_i (\bar x ,  y'(t)\cdot e_i,t) \leq C h^2.
\end{equation*}

\medskip

This implies that the value function $u$ is semi-concave with respect to~$x$ with a linear modulus of semi-concavity,  locally in $J_i\setminus \{\vertex\}$ and for $t$ bounded away from the horizon $T$. More precisely,
consider $t\in[0,T)$ and $x,y\in J_i\setminus\{\vertex\}$ with $0<r\leq |x|,|y|\leq R$.
there exists a constant $C$ (depending in particular on $r$, $R$ and $T-t$) such that
$\lambda u(x,t)+(1-\lambda)u(y,t)-u(\lambda x+(1-\lambda) y,t)\leq C\lambda(1-\lambda)|x-y|^2$,  for every $\lambda \in[0,1]$.

In the remainder of the paper, the assumption [H+] will not be made because we will not need the local semi-concavity of $u$ with a linear modulus.
\end{remark}

\subsubsection{Regularity of $u$ along optimal trajectories}\label{sec:regularity_u}
Here, we investigate some regularity properties of $u$ in the interiors of the edges.\\
For $t\in[0,T)$, $x\in J_i\setminus \{\vertex\}$  for some $i\in \{1,\dots,N\}$, let us set $$\mathcal{A}(x,t)=\left\{z'(t)\cdot e_i\;:\; z\in \Gamma^{\rm{opt}}_t[x] \right\}$$ (recall that $z'(t)$ is well defined for all $ z\in \Gamma^{\rm{opt}}_t[x]$ because $x\not= \vertex$).% see point 4 in Lemma \ref{lemma:EL}).

\begin{lemma}\label{lemma:note4.9}
Consider $t\in[0,T)$, $x\in J_i\setminus \{\vertex\}$,  $y\in \Gamma^{\rm{opt}}_t[x]$, and  $t_*=T\wedge \min \{\tau\in[t,T]\;:\; y(\tau)=\vertex\}$. Set $x=\bar x e_i$ and $y(s)=\bar y (s) e_i$ for $s\in (t,t_*)$.

Under assumptions [H1] to [H5], % \ref{ass_for_section_2_4} and \ref{ass_for_section_2_7},
the following properties hold:
\begin{itemize}
\item[(i)] For every $s\in(t,t_*)$, $\mathcal{A}(y(s),s)$ is singleton.

\item[(ii)] $D_x u(x,t)$ exists if and only if $\mathcal{A}(x,t)$ is a singleton. In this case,
\begin{equation}
\label{eq:10003}
%  \begin{split}
\mathcal{A}(x,t) =   (\partial_a\ell_i(\bar x,\cdot, t)  )^{-1}  \left(  \Bigl\{  -D_x u(x,t)  \Bigr\} \right)   
=    \Bigl\{   - \partial_p H_{i}\left( \bar x, D_x u(x,t),   t\right)    \Bigr\}.
% \end{split}
\end{equation}
\item[(iii)] %For any   $y\in \Gamma^{\rm{opt}}_t[x]$, set $t_*=T\wedge \min \{\tau\in[t,T]\;:\; y(\tau)=\vertex\}$.
For all $s\in(t,t_*)$, $u(\cdot,s)$ is differentiable at $y(s)$  with
$D_x u(y(s),s) = -\partial_a \ell_i (\bar y(s) ,   \bar y'(s),s)$. 
\end{itemize}
\end{lemma}
\begin{proof}
{\sl(i)}  The arguments are similar to those in the proof of \cite[Lemma 4.9-(1)]{C}.
%, so we refer the reader to that paper for the details.
For every $s\in(t,t_*)$,  $z\in \Gamma^{\rm{opt}}_s[y(s)]$, let $y_h$ be the admissible trajectory defined as follows:
\begin{equation*}
y_h(\tau)=\left\{\begin{array}{ll}
y(\tau)&\qquad\textrm{if }\tau \in[t,s-h],\\
y(s-h)+ \frac{\tau -s +h}{2h} (z(s+h)-y(s-h)) &\qquad\textrm{if }\tau \in(s-h, s+h),\\
z(\tau)&\qquad\textrm{if }\tau \in[s+h,T].
\end{array}\right.
\end{equation*}
The desired result is obtained by comparing the costs associated with $y_h$ and with the concatenation of $y$ and $z$ at time $s$, then by letting $h\to 0$.

{\sl (ii)} Assume that $D_x u(x,t)$ exists. We wish to prove that $\mathcal{A}(x,t)$ is a singleton.
For $z\in \Gamma^{\rm{opt}}_t[x]$, let us set $ z'(s) =\bar\alpha(s)  e_i$ for $s\in [t,t_*)$ with
$t_*=T\wedge \min \{\tau\in[t,T]\;:\; z(\tau)=\vertex\}$.
From Lemma \ref{lemma:CC2_3.1}, $\ds -\partial_a \ell_i (\bar x , \bar\alpha(t),t)\in D^+_{x}u(x,t)$ where 
$D^+_{x}u(x,t)$ is defined in \eqref{eq:super_diff_u_x}.
%= \{q\in \R: \limsup_{J_i\ni y\to x}\frac{u(y,t)-u(x,t)-q  ( y- x)\cdot e_i}{| y- %x|}\leq 0\}  $.
Moreover, since $u(\cdot, t)$ is differentiable at $x$, $D^+_x u(x,t)$ is a singleton. Hence, $\{ \partial_a \ell_i (\bar x ,  e_i\cdot z'(t),t), z\in  \Gamma^{\rm{opt}}_t[x]\}$
is the singleton $\{-D_x u(x,t)\}$. Since $\ell_i$ is strictly convex with respect to its second argument, this implies that $\mathcal{A}(x,t)$ is the singleton
defined in the first identity in \eqref{eq:10003}. The second identity in \eqref{eq:10003} comes from the definition of $H_i$.
 
\medskip

Conversely,  assume that $\mathcal{A}(x,t)$ is a singleton. 
%In order to obtain the differentiability of $u$  w.r.t. $x$ at $(x,t)$, 
Let us prove that
\begin{displaymath}
D^*u(x,t)=\left\{
p\in \R:  \exists (x_n)_{n\in \N}  \hbox{ with} \left| \begin{array}[c]{ll} & \ds \lim_{n\to\infty} x_n= x, \\
 &u(\cdot, t) \hbox{ is differentiable at } x_n ,    
\\ & \ds \lim_{n\to \infty} D_x u(x_n,t)  =p
\end{array}\right.
\right \}  
\end{displaymath}
is a singleton. Since $u$ is semi-concave from the last point in Lemma \ref{lemma:CC2_3.1},
%Proposition \ref{prop:semiconc_loc},
this will imply that $u(\cdot, t)$ is differentiable w.r.t. $x$ at $(x,t)$, from 
\cite[Theorem 3.3.6, eq.(3.20)]{MR2041617} and
from \cite[Proposition 3.3.4(d)]{MR2041617}.

Consider $p\in D^*u(x,t)$ and a sequence $(x_n)_n$, $x_n\in J_i\setminus\{\vertex\}$ as in the definition of $D^*u(x,t)$.
There exist $y_n\in \Gamma^{\rm{opt}}_t[x_n]$ and some  $\delta >0 $ and  $\bar t>t $ such that, for all $n$, $y_n(s)\in J_i\setminus\{\vertex\}$ 
and $e_i\cdot y_n(s)>\delta$ for all $s\in [t,\bar t]$. Set $x_n=\bar x_n e_i$, $ y_n(s)=\bar y_n(s)  e_i$
and $\bar\alpha_n= \bar y_n'(s)$
for $s\in [t,\bar t]$. 

From the necessary optimality condition in Lemma \ref{lemma:EL} and the differentiability of $u(\cdot, t)$ at  $x_n $,  we know that
\begin{equation} \label{eq:10004}
\begin{split}
\frac d {ds} \partial_a \ell_i (\bar y_n(s),\bar\alpha_n(s),s) =  \partial_x \ell_i (\bar y_n(s),\bar\alpha_n(s),s),&\quad s\in [t,\bar t],\\
\bar y_n(t)=  \bar x_n,&\\
\bar\alpha_n(t)= -\partial_p H_i(\bar x_n,   D_x u(x_n,t) , t).
\end{split}
\end{equation}
From the bounds on optimal trajectories, we deduce that, up to the extraction of a subsequence, $y_n$  converges uniformly to some function $z$ as $n\to\infty$, and weakly in $W^{1,p}(t,T;\R^N)$,
and that $z$ belongs to $ \Gamma^{\rm{opt}}_t[x]$. Since $z(s)\in  J_i\setminus\{\vertex\}$ and $e_i\cdot z(s)\ge \delta$  for all $s\in [t,\bar t]$, we can set $ z(s)=\bar z(s) e_i$ and
$ \bar z'(s)=\bar\alpha(s)$  for all $s\in [t,\bar t]$. We obtain
\begin{equation*}
\begin{split}
\frac d {ds} \partial_a \ell_i (\bar z(s),\bar\alpha(s),s) =  \partial_x \ell_i (\bar z(s),\bar\alpha(s),s),&\quad s\in [t,\bar t],\\
   \bar z(t)= \bar x, \end{split}
\end{equation*}
because $z\in \Gamma^{\rm{opt}}_t[x]$.
We wish to prove that $\partial_a \ell_i(\bar x, \bar\alpha(t), t)=-p$.
This will characterize $p$ because ${\mathcal {A}}(x,t)$  a singleton,
 hence show $u(\cdot, t)$ is differentiable at $x$ with $D_x u(x,t)=-\partial_a \ell_i(\bar x,  \bar z'(t), t)$ for every $z\in\Gamma^{\rm{opt}}_t[x]$.

From Theorem~\ref{optcur_lip},  there exists $V>0$ such that $\|\alpha_n\|_{L^\infty(t,\bar t)} \le V$  for all $n$, and $\|\alpha\|_{L^\infty(t,\bar t)} \le V$.
Hence, $ (\partial_x \ell_i (\bar y_n(\cdot),\bar\alpha_n(\cdot),\cdot))_n$ is bounded in $L^\infty(t,\bar t)$ and, from \eqref{eq:10004}, $ (\partial_a \ell_i (\bar y_n(\cdot),\bar\alpha_n(\cdot),\cdot))_n$ is bounded in $W^{1,\infty}(t,\bar t)$. 
Possibly after a further extraction of a subsequence, 
$ \partial_a \ell_i (\bar y_n(s),\bar\alpha_n(s),s)$ converges to some function $\pi$  weakly-$*$  in $W^{1,\infty} (t,\bar t)$ and in $C([t,\bar t] )$.
But also $ \partial_a \ell_i (\bar z(s),\bar\alpha_n(s),s)$ converges to $\pi$  in $C([t,\bar t] )$,  because $ \partial_a \ell_i (\zeta,\cdot,s)$ 
depends continuously on $(\zeta, s)\in \R_+\times [0,T]$ and $\bar y_n$ converges uniformly to $\bar z$ in $[t,\bar t]$.  Hence, $\bar\alpha_n$
converges to $s\mapsto  (\partial_a \ell_i (\bar z(s),\cdot,s))^{-1}(\pi(s))$ uniformly on $[t,\bar t]$. This implies that $\bar\alpha(s)= (\partial_a \ell_i (\bar z(s),\cdot,s))^{-1}(\pi(s))$ and that $\bar\alpha_n$ converges to $\bar\alpha$ uniformly on $[t,\bar t]$.
From this and \eqref{eq:10004}, we deduce that 
$p=-\partial_a \ell_i(\bar x, \bar\alpha(t), t)$, or equivalently that
$\bar\alpha(t)= -\partial_p H_i(\bar x, p, t)$.
Hence, being ${\mathcal {A}}(x,t)$  a singleton,  $Du^* (x,t)$ is  also a singleton. The proof of (ii) is complete.

(iii) is obtained by combining the first two points.
\end{proof}

\begin{remark}\label{rmk:h7} Note that the assumptions [H1] to [H5] are not enough to guarantee the uniqueness of local solutions to the  Cauchy problem originating from the optimality conditions, namely:
given $\bar x>0$ and $\bar \alpha\in \R$,
find $\bar y$ such that
\begin{equation*}
\label{eq:opt_cond1_2}
\begin{split}
\frac d {ds} \left(\partial_a \ell_i(\bar y(s),\bar y'(s),s)\right)=\partial_y \ell_i(\bar y(s),\bar y '(s),s), \quad \hbox {for } s>t,\\
\bar y (t)=\bar x, \quad \hbox{and} \quad    \bar y'(t)=\bar \alpha.
\end{split}
\end{equation*}
In order to obtain  uniqueness of local solutions,  further assumptions are needed, for example
\begin{description}
\item{[H6]}
\begin{enumerate}
\item For  every $s\in [0,T]$,  $\partial^2 _{a,a} \ell_i(\cdot,\cdot,s)$,  $\partial^2 _{x,a} \ell_i(\cdot,\cdot,s)$, $\partial^2 _{s,a} \ell_i(\cdot,\cdot,s)$  are continuous functions defined on $\R_+\times \R$ 
\item   There exists a constant $c>0$ such that $\partial^2 _{a,a} \ell_i  \ge c$
\item  For every compact $K $ of $\R_+\times \R$, there exists a function $\Lambda_K\in L^1(0,T)$ such that for all $(x,a)\in K$, $(\tilde x, \tilde a) \in K$ and for almost all $s\in [0,T]$,
\begin{equation*}
\begin{split}
|\partial^2 _{a,a} \ell_i (x,a,s) - \partial^2 _{a,a} \ell_i (\tilde x,\tilde a,s)| \le \Lambda_K(s) (  |x-\tilde x|+|a-\tilde a |),\\
|\partial^2 _{x,a} \ell_i (x,a,s) - \partial^2 _{x,a} \ell_i (\tilde x,\tilde a,s)| \le \Lambda_K(s) (  |x-\tilde x|+|a-\tilde a |),\\
|\partial^2 _{s,a} \ell_i (x,a,s) - \partial^2 _{s,a} \ell_i (\tilde x,\tilde a,s)| \le \Lambda_K(s) (  |x-\tilde x|+|a-\tilde a |).
\end{split}
\end{equation*}
\end{enumerate}
\end{description}
If assumptions [H1] to [H6] %[H1] to [H7]
hold, then stronger results may be obtained, see 
\cite[proofs of Lemma 2.27, Corollary 2.28, Lemma 2.30]{AMMT3}:
\begin{enumerate}
\item Point (i) in Lemma \ref{lemma:note4.9} is improved as follows:
for every $s\in(t,t_*)$, $y'|_{(s,t_*)}$ is the unique optimal control for $u(y(s),s)$ up to time $t_*$. 
\item  Consider two trajectories $\gamma_1\in\Gamma^{\textrm{opt}}[x_1]$ and $\gamma_2\in\Gamma^{\textrm{opt}}[x_2]$ such that $\gamma_1(t)=\gamma_2(t)\in J_i\setminus\{\vertex\}$ for some $t\in (0,T)$. Let $I_k$,  $k=1,2$, be the largest open interval containing $t$ such that $\gamma_k(s)\in J_i\setminus\{\vertex\}$ for $s\in I_k$. Then $I_1=I_2$.
\item (Optimal synthesis) Consider $t\in[0,T)$, $x\in J_i\setminus \{\vertex\}$ for some $i\in \{1,\dots,N\}$. 
 \begin{enumerate}
 \item  \label{pointa} If there exist $t_*\in (t,T]$ and $\bar y\in W^{1,\infty}(t,t_*)$ such that $\bar y(t)=|x|$, $\bar y(s)>0$ in $[t,t_*)$,
 and $\bar y'(s)=- \partial_p H_i (\bar y(s), D_x u(\bar y(s) e_i,s) , s) $ a.e. in $(t,t_*)$,
 then the control 
 \begin{equation*}
\alpha(s)= - \partial_p H_i (\bar y(s), D_x u(y(s),s) ,s) , s) e_i
\end{equation*}
is optimal for $u(x,t)$ in $[t,t_*]$.
\item If $u(\cdot,t)$ is differentiable at $x$, then there exist a unique $t_*$ and $\bar y$ as in point \ref{pointa} and such that $t_*=T$ or $\bar y (t_*)=0$.
\end{enumerate}
\end{enumerate}
\end{remark}

\section{ Mean Field Games}\label{sect:MFGequil}
%
%Setting and notations
%
We are now ready to tackle mean field games.
Relying on the results contained in Section~\ref{sec:OC}, we investigate the existence of relaxed equilibria and study their properties.

\subsection{Setting and notations}\label{subsec:setting}
In what follows, we introduce some notation and the notion of relaxed MFG equilibria  (in Definition  \ref{sec:setting-notation-3} below). The latter are not described by a system of PDEs but rather by probability measures on a set of admissible trajectories satisfying a suitable condition (see \eqref{eq:43}), given the initial distribution of states $m_0$. The relevant  sets of trajectories, denoted $\Gamma_C$ 
and $\Gamma^{\textrm{Lip}}_{C,V}$ defined below, respectively contain trajectories  uniformly bounded in $W^{1,p}(0,T; \R^N) \cap C([0,T]; \cG) $ and trajectories  uniformly bounded in ${\rm Lip} ([0,T]; \cG)$. 
The notion of relaxed equilibria is described in Definition  \ref{sec:setting-notation-3} below, and all the necessary ingredients  are defined in the following paragraphs.

\paragraph{Probability sets}
Hereafter, $m_0$ is a fixed probability measure on $\cG$ with compact support.

Given a topological space $X$, let $\cP(X)$ be the set of Borel probability measures on~$X$.

Let $\cP_1(\cG)$ denote the set of Borel probability measures on $\cG$ with finite first moment,
endowed with the Kantorovich-Rubinstein distance $d_1$. 

Recall that for some fixed $p\in (1,+\infty)$, the set of trajectories $\Gamma$ is defined in Subsection \ref{setting}. We endow $\Gamma$  with the topology of uniform convergence in $[0,T]$. %The set $\cP(\Gamma)$ will play an important role in the following.
For every $\mu\in \cP(\Gamma)$ and  $t\in [0,T]$,  the marginal $m^\mu(t)\in \cP( \cG)$ is defined by $m^\mu(t)=e_t\sharp \mu$,  where  $e_t:\Gamma\to \cG$ is  the evaluation map defined by $e_t(y)=y(t)$.

\paragraph{Costs} The running cost and the terminal cost depend on the distribution of states.
For $i=1,\dots, N$, let us introduce  $L_i\in C(\cP_1(\cG); C^1( \R_+\times\R))$, $\cP_1(\cG) \ni m\mapsto L_i[m]\in  C^1( \R_+\times\R)$: if the sequence $(m_n)_n $,  $m_n\in \cP_1(\cG)$, and $m\in\cP_1(\cG)$ are such that $\lim_{n\to\infty} d_1(m_n,m)=0$, then for each compact subset $K$ of $\R_+\times \R$, 
$L_i[m_n]|_{K}$ converges to $L_i[m]|_{K}$  in $C^1 (K)$. Similarly, we consider $G_i\in C_b(\cP_1(\cG); C_b( \R_+))$, $L_* \in   C_b(\cP_1(\cG); \R)$  and $G_*\in  C_b(\cP_1(\cG); \R)$.
Hereafter, we will always assume that $L_i$ and $L_*$ satisfy the following minimal set of assumptions (the counterparts of assumptions [H1] to [H4] in Section \ref{sec:OC}):
\begin{itemize}
\item[[$H_{MFG,1}$]]
\begin{enumerate}
\item There exists a positive constant $C_0$ such that for all $i\in \{1,\dots,N\}$, $m\in \cP_1(\cG)$, $x\in \R_+$, $a\in \R$, $L_i[m]\ge C_0(|a|^p-1)$ and  $L_* [m]\ge -C_0$
\item For any $r>0$, there exists a modulus $\omega_r$ such that
for all $i\in \{1,\dots,N\}$, $m_1,m_2\in \cP_1(\cG)$, $x\in [0,r]$, $a\in \R$, 
$  | L_i[m_1](x,a)- L_i[m_2](x,a)|\le (1+|a|^p) \omega_r(d_1(m_1,m_2))$
\item For any $r>0$, there exists $\tilde c (r)\ge 0$ such that
for all $i\in \{1,\dots,N\}$, $m\in \cP_1(\cG)$, $x\in [0,r]$, $a\in \R$,
$|\partial_x L_i[m](x,a)|\le \tilde c (r) (|a|^p+1)$ and  $|\partial_a L_i[m](x,a)|\le \tilde c (r) (|a|^{p-1}+1)$
\item For all  $m\in \cP_1(\cG)$, $x\in [0,r]$, $a\mapsto L_i[m](x,a)$ is strictly convex.
\end{enumerate}
\end{itemize}
We will sometimes make the following assumption (the counterpart of assumption [H5] in Section \ref{sec:OC}):
%We will sometimes make the following assumptions (the counterparts of assumptions [H5] and [H6] in Section \ref{sec:OC}):
\begin{itemize}
\item[[$H_{MFG,2}$]] For all $i\in \{1,\dots, N\}$, $m\in \cP_1(\cG)$, $G_i[m] \in C^1( \R_+)$.
\end{itemize}
  For $m\in \cP_1(\cG)$, a trajectory $ y\in \Gamma$, let us set   $\bar y =|y|$  and define
\begin{eqnarray*}
L_\vertex[m]=\min\left\{L_*[m],\min_{i=1,\dots,N}L_i[m](0,0)\right\}, \notag \\
L[m](y(t),y'(t)) = \sum_{i=1}^N L_i[m](\bar y(t) ,\bar y'(t))\car_{y(t)\in J_i\setminus\{\vertex\}} + L_\vertex[m]\car_{y(t)=\vertex}, \label{eq:37b}\\
G[m](x)=\ds \sum_{i=1}^N G_i[m](|x|)\car_{x\in J_i\setminus\{\vertex\}} +\min\left\{G_*[m],\min_{i=1,\dots,N}G_i[m](0)\right\}\car_{x=\vertex}.
\end{eqnarray*}

\paragraph{Example of running costs}
Consider
\begin{itemize}
\item an exponent $p>1$
\item for  $i=1,\dots,N$ and $j=1,2$, a function $\kappa_{i,j}\in C^\infty_0( [0,+\infty) )$  such that  $0 \le  \kappa_{i,j}\le 1$ and
  $\kappa_{i,j}'(0)  =0$.
%  and  $\kappa_{i,1}$ is non identically $0$
\item  for $i=1,\dots,N$ and $j=1,2$, a smooth and bounded function $h_{i,j}$ defined on $J_i\times [0,1]$,
$h_{i,1}$ being bounded from below by a positive constant.
\end{itemize}
Then, for $i=1,\dots,N$ and $j=1,2$,  $y=\bar y e_i$, $a\in \R$, set
\begin{displaymath}
  K_{i,j}[m](y)=  h_{i,j} \left(    y,   \int_{\cG}   \kappa_{i,j}( d(y,z))m(dz)\right)  
\end{displaymath}
  and define
\begin{displaymath}
L_i[m](\bar y ,a)= K_{i,1}[m](y) |a|^p +K_{i,2}[m](y).
\end{displaymath}
Similarly, set $L_*[m]= h_{*} \left(   \int_{\cG}   \kappa_{*}( d(0,z))m(dz)\right)$, where $h_* \in C^\infty([0,1])$, and $\kappa_{*}\in C^\infty_0( [0,+\infty) )$ is such that $0 \le  \kappa_{*}\le 1$.
\\
The functions $L_*$ and $L_i$, $i=1,\dots, N$ fulfill the above assumptions.

\paragraph{Admissible trajectories}
For $C\ge 0$, let us introduce the set
\begin{equation}\label{eq:gamma_tilde_x}
  \Gamma_C[x]=\left\{y\in \Gamma[x]
\;:\; \|y\|_{L^\infty(0,T)}\leq C,\  \| y'\|_{L^p(0,T)}\leq C\right\}
\end{equation}
which is not empty if $C\ge |x|$,
and $\Gamma_C=\bigcup_{x\in\cG} \Gamma_C[x]$, endowed with the topology of uniform convergence on $[0,T]$. For every positive constant $C$, $\Gamma_C$ is compact.

\paragraph{Lipschitz admissible trajectories}
Given two positive constants $V$ and $C$, set
\begin{equation*}\label{gamma_lip}
\Gamma^{\textrm{Lip}}_{C,V}[x]=\left\{y\in \Gamma_C[x]\;:\; \|y'\|_{L^\infty(0,T)}\leq V\right\},\quad \hbox{and} \quad 
\Gamma^{\textrm{Lip}}_{C,V}=\bigcup_{x\in\cG}\Gamma^{\textrm{Lip}}_{C,V}[x],
\end{equation*}
and endow $\Gamma^{\textrm{Lip}}_{C,V}$ with the topology of uniform convergence. The set $\Gamma^{\textrm{Lip}}_{C,V}$ is compact.

\paragraph{The set $\cP( \Gamma_C)$ and the associated costs}
The set $\cP(\Gamma_C)$ of probability measures on $\Gamma_C$ is endowed with the weak topology.
 If $\mu\in \cP(\Gamma_C)$, then $\textrm{supp}(m^\mu(t))\subset\{x\in\cG\;:\;|x|\leq C\}$, and the map $t\mapsto m^\mu(t)$ belongs to $C^{(p-1)/p}([0,T];\cP_1(\cG))$.
Hence, for $(x,a)\in\R_+\times\R$, the function $t\mapsto L_i[m^\mu(t)](x,a)$ is continuous and  then $(t,x,a) \mapsto  L_i[m^\mu(t)](x,a)$ satisfies the assumptions [H1] to [H4] because the functions $L_i$ satisfy the assumption [$H_{ MFG,1}$].

With $\mu\in \cP_1(\Gamma)$ and $y\in\Gamma[x]$, we associate the cost
\begin{equation}\label{costMFG}
J^\mu(x;y)
=\int_0^T L[m^\mu(\tau)] (y(\tau),y'(\tau))d\tau+  G[ m^\mu(T)](y(T)).
\end{equation}

\paragraph{Optimal trajectories} 
Let $C$ and $\widetilde C$ be positive numbers. Fix $\mu\in \cP(\Gamma_{\widetilde C})$. For every $x\in \cG$, let us set 
\begin{equation}\label{eq:42new}
\GammamuoptC[x]=\left\{y\in \Gamma_C[x]\;:\; J^\mu(x;y)=\min_{ \widetilde y\in \Gamma[x]}  J^\mu(x; \widetilde y) \right\},
\end{equation}
where $J^\mu $ is defined in (\ref{costMFG}). Proposition~\ref{prp:ex_OT} implies that for each $\mu\in \cP(\Gamma_{\widetilde C})$ and $x\in \cG$,  $\GammamuoptC[x]$  is not empty if $C$ is sufficiently large (depending on $\widetilde C$).

\paragraph{The sets  $\cP_{m_0}(\Gamma_C)$ and $\cPmzeroGammaLipCV$} Let $\cP_{m_0}(\Gamma_C)$ (resp.  $\cPmzeroGammaLipCV$)  denote the set of measures $\mu\in\cP(\Gamma_C)$
(resp. $\mu\in\cP( \Gamma^{\mathrm{Lip}}_{C,V})$)
such that $e_0\sharp \mu=m_0$. For $C$ and $V$ sufficiently large,
 $\cP_{m_0}( \Gamma_C)$
and $\cPmzeroGammaLipCV$ are not empty.

\paragraph{Relaxed equilibria}

\begin{definition}
  \label{sec:setting-notation-3}
For $C>0$, the complete probability measure $\mu \in \cP_{m_0}(\Gamma_C)$ is a relaxed equilibrium associated with the initial distribution $m_0$ if 
\begin{equation}\label{eq:43}
{\rm{supp}}(\mu)\subset \mathop \bigcup_{
x\in {\rm{supp}}(m_0)}\GammamuoptC[x].
\end{equation}
\end{definition}

\subsection{Existence of a relaxed equilibrium}
The property of $(x,\mu)\rightrightarrows \Gamma^{\mu,\rm{opt}}_C[x]$ obtained in the following proposition is important for proving the existence of relaxed equilibria because it allows one to apply Kakutani fixed point theorem, see, for example, \cite{AMMT3, CC}.
\begin{proposition}\label{prop2}
We make the assumption [$H_{ MFG,1}$].
Let $C\geq 0$ be such that $\cP_{m_0}(\Gamma_C)\not= \emptyset$. Consider $\mu\in\cP_{m_0}(\Gamma_C)$, $x\in\textrm{supp}(m_0)$,  a sequence $(\mu_n)_{n\in\N}$, $\mu_n\in\cP_{m_0}(\Gamma_C)$,  weakly convergent to $\mu$ as $n\to\infty$ and a sequence of points $(x_n)_{n\in\N}$,  $x_n\in\cG$ converging to $x$ as $n\to\infty$. 
If $(y_n)_{n\in\N}$, $y_n\in \Gamma_C^{\mu_n,{\rm opt}}[x_n]$  converges uniformly to some $y$ as $n\to\infty$, then $y$ belongs to $\GammamuoptC[x]$.
\end{proposition}
\begin{proof}
Recall that $\lim_{n\to \infty} \max_{t\in [0,T] }  d_1(m^{\mu_n}(t), m^\mu(t)) =0$, and that since $\Gamma_C[x]$ is closed,   $y\in \Gamma_C[x]$, and the sequence  $(y'_n)_{n\in\N}$  converges weakly to $y'$ in $L^p(0,T;\R^N)$.
We wish to prove that $J^\mu(x;y)\leq J^\mu(x;\hat y)$ for every $ \hat y \in\Gamma[x]$. Fix $\hat y\in\Gamma[x]$ and consider a sequence $(\hat y_n)_n$, $\hat y_n\in\Gamma[x_n]$ as in Lemma~\ref{lemma1}. Since $y_n\in \Gamma_C^{\mu_n,{\rm opt}}[x_n]$, 
\begin{equation}\label{eq:34}
J^{\mu_n}(x_n;y_n)\leq J^{\mu_n}(x_n;\hat y_n).
\end{equation}
Let us focus on the left-hand side of~\eqref{eq:34}.  First,
$|G[m^{\mu_n}(T)](y_n(T))-G[m^{\mu}(T)](y(T))|$ is not larger than
$
|G[m^{\mu_n}(T)](y_n(T))-G[m^{\mu_n}(T)](y(T))|+|G[m^{\mu_n}(T)](y(T))-G[m^{\mu}(T)](y(T))|$, hence vanishes as $ n\to \infty$.
On the other hand,
\begin{equation*}
  \int_0^TL[m^{\mu_n(\tau)}](y_n(\tau),y_n'(\tau))d\tau-\int_0^TL[m^{\mu(\tau)}](y(\tau),y'(\tau))d\tau
= A_n + B_n,
\end{equation*}
where
\begin{eqnarray*}
A_n &=&\int_0^TL[m^{\mu_n(\tau)}](y_n(\tau),y_n'(\tau))d\tau-\int_0^TL[m^{\mu(\tau)}](y_n(\tau),y_n'(\tau))d\tau,\\
B_n &=&\int_0^TL[m^{\mu(\tau)}](y_n(\tau),y_n'(\tau))d\tau-\int_0^TL[m^{\mu(\tau)}](y(\tau),y'(\tau))d\tau.
\end{eqnarray*}          
From the assumptions on the costs, in particular point 2. in [$H_{ MFG,1}$], we see that $\lim_{n\to \infty} A_n=0$. On the other hand, Lemma~\ref{lem:LSC} implies $\liminf_n B_n\geq 0$. Combining  the observations above yields
\begin{equation}\label{eq:34quat}
\liminf_{n\to\infty}  J^{\mu_n}(x_n;y_n)\geq J^\mu(x;y).
\end{equation}

Consider now the right-hand side of~\eqref{eq:34}.
Using the assumptions on the costs, arguments similar to those in the proof of Lemma~\ref{lemma1} lead to:
\begin{equation}\label{eq:sum1}
\lim_{n\to \infty}  J^{\mu_n}(x_n;\hat y_n) = J^{\mu}(x;\hat y).
\end{equation}
In conclusion, \eqref{eq:34}, \eqref{eq:sum1} and~\eqref{eq:34quat} yield
$
J^\mu(x;y)\leq J^\mu(x;\hat y)$.
Since $\hat y$ is
arbitrary, we get
$J^\mu(x;y)=\min_{\hat y\in \Gamma[x]}J^\mu(x;\hat y)$.
\end{proof}

We can now state our main result about the existence of relaxed equilibria.
\begin{theorem}\label{sec:setting-notation-6}
\begin{itemize}
\item[($i$)] Under the assumption [$H_{ MFG,1}$], for $C$ large enough, there exists a relaxed equilibrium $\mu\in \cP_{m_0}(\Gamma_C)$
\item[($ii$)] Under the assumptions [$H_{ MFG,1}$] and [$H_{ MFG,2}$], for $C$ and $V$ large enough, there exists a relaxed equilibrium $\mu\in \cPmzeroGammaLipCV$. 
\end{itemize}
\end{theorem}
\begin{proof}
We  describe the main steps in the proof and refer the reader to \cite{AMMT3} for more details, because point $(i)$ (resp. point $(ii)$) relies on the same arguments as  \cite[Proof of Theorem 3.5]{AMMT3} (resp. \cite[Proof of Theorem 3.6]{AMMT3}) with easy adaptations, replacing \cite[Proposition 3.10]{AMMT3} with Proposition~\ref{prop2}.\\
$(i)$. Let the multivalued map $E:\cP_{m_0}(\Gamma_C)\rightrightarrows \cP_{m_0}(\Gamma_C)$ be defined as follows:
\begin{equation}\label{eq:mappa}
E(\mu)=\{\hat \mu\in\cP_{m_0}(\Gamma_C)\;:\; \textrm{supp}\,\hat \mu_x\subset \GammamuoptC[x]\quad m_0-\textrm{a.e. } x\in\cG\},
\end{equation}
where $\{\hat \mu_x\}_{x\in\cG}$ is the family of Borel probability measures on $\cP_{m_0}(\Gamma_C)$ obtained applying the disintegration theorem \cite[Theorem 5.3.1]{AGS} in order to disintegrate $\hat \mu$  with respect to $m_0$. The proof consists of checking that $E$ admits a fixed point, thanks to Kakutani theorem. To this end, it suffices to prove the following properties
\begin{itemize}
\item[(a)] for every $\mu\in\cP_{m_0}(\Gamma_C)$, the set $E(\mu)$ is not empty and convex
\item[(b)] the map~$E$ has a closed graph.
\end{itemize}
The non-emptiness is proved using the closed graph property in Proposition~\ref{prop2} and a Borel selection that associates to each point~$x\in\cG$ an optimal trajectory. Convexity is proved using the disintegration theorem \cite[Theorem 5.3.1]{AGS}. On the other hand, the closed graph property of the map~$E$ is proved using Kuratowski theorem (\cite[Proposition 5.1.8]{AGS}) and again Proposition~\ref{prop2} and the disintegration theorem \cite[Theorem 5.3.1]{AGS}.\\
$(ii)$. This proof follows the same steps as in  point~$(i)$. 
Consider the restriction $E_0$ of the  multivalued map~$E$ to the set $\cPmzeroGammaLipCV$ where $C$ and $V$ are chosen sufficiently large. 
To prove that $E_0$ also fulfills the hypotheses of Kakutani fixed point theorem,  we need to check that
\begin{itemize}
\item[(a)] $E_0(\mu)\subset \cPmzeroGammaLipCV$, $\forall\mu\in \cPmzeroGammaLipCV$
\item[(b)] $\cPmzeroGammaLipCV$ is compact
\item[(c)] $E_0(\mu)$ is a non-empty convex set, $\forall\mu\in \cPmzeroGammaLipCV$
\item[(d)] $E_0$ has the closed graph property.
\end{itemize}
Point~$(a)$ relies on Theorem~\ref{optcur_lip}. Point~$(b)$ follows from the compactness of~$\Gamma^{\textrm{Lip}}_{C,V}$ and Prokhorov theorem \cite[Theorem 5.1.3]{AGS}. 
The two points~$(c)$ and~$(d)$ are proved exactly as (a) and (b)  in point~$(i)$.
\end{proof}
\begin{remark}
Under assumptions [$H_{ MFG,1}$] and [$H_{ MFG,2}$], Proposition~\ref{prp:ex_OT} and Theorem~\ref{optcur_lip}  ensure that each relaxed equilibrium $\mu$ belongs to $\cPmzeroGammaLipCV$, for $C$ and $V$ large enough. 
\end{remark}
\subsection{Mild solutions}
\begin{definition}\label{def:mildsol}
Let $\mu\in\cP_{m_0}(\Gamma_C)$ be a relaxed equilibrium.
The pair $(u,m)$ is the associated mild solution if
$m\in C([0,T];\cP_1(\cG))$ is defined by $m(t)=m^\mu(t)=e_t\#\mu$ and 
$u$ is the value function  associated with $\mu$ by \begin{equation}\label{eq:4.1}
u(x,t)= \inf_{ y\in \Gamma_t[x]} 
J^\mu_t(x; y ),
\end{equation}
where $J^\mu_t(x;y)=\int_t^T L[m(\tau)] (y(\tau),y'(\tau))d\tau+  G[ m(T)](y(T))$.
\end{definition}
\begin{remark}\label{rmk:reg_m}
For $\mu\in\cP_{m_0}(\Gamma_C)$ (resp. $\mu\in \cPmzeroGammaLipCV$), using 
Proposition \ref{prp:ex_OT},
it can be checked as in  \cite[Lemma 3.8]{AMMT3}  that the function $t\mapsto m(t)$ belongs to $C^{(p-1)/p}([0,T];\cP_1(\cG))$ (resp. ${\rm Lip}([0,T];\cP_1(\cG))$).
As a consequence, the functions $L_i[m(\cdot)](x,a)$ belong to $C^{(p-1)/p}([0,T])$ (resp. to $Lip([0,T])$).
\end{remark}
Proposition \ref{mild_sol_u} and Lemma \ref{lemma:reg_u_mfg} that follow contain direct consequences of the results proved in Section \ref{sec:OC}.
\begin{proposition} \label{mild_sol_u}
Let us make the assumption [$H_{ MFG,1}$]. Let $\mu\in\cP_{m_0}(\Gamma_C)$ be a relaxed equilibrium and $(u,m)$ be the related mild solution. There holds:
\begin{enumerate}
\item (Dynamic programming principle) For all $0\le t\le \bar t \le T$, $x\in \cG$,
\[
u(x,t)= \inf_{y\in \Gamma_{t,\bar t}[x]}\left\{
u(y(\bar t),\bar t)+\int_t^{\bar t} L[m(\tau)](y(\tau),y'(\tau))d\tau
\right\}
\]
\item The function $u$ is continuous in~$\cG\times[0,T)$. For all $x\in \cG$, $t\mapsto u(x,t)$ is continuous on $[0,T]$, and $u(x,T)=G[m(T)](x)$
\item $u$  is a viscosity solution to \eqref{HJ} where the Hamiltonians $H_i$ and $H_\vertex$ are the ones given at the beginning of Section \ref{SHJB} replacing $\ell_i(\cdot,\cdot,t)$ and $\ell_0(t)$ respectively with $L_i[m(t)]$ and $L_0[m(t)]$.  
\end{enumerate}
Under the assumptions [$H_{ MFG,1}$] and [$H_{ MFG,2}$],
\begin{description}
\item{4.a.} $u$ is locally Lipschitz continuous in~$\cG\times[0,T)$
\item{4.b.} If moreover,  $G$ maps $\cP_1(\cG)$   to $C (\cG)$, then $u$ is locally Lipschitz continuous in~$\cG\times[0,T]$.
\end{description}
\end{proposition}

\begin{lemma}\label{lemma:reg_u_mfg}
  Let us make the assumptions [$H_{ MFG,1}$] and [$H_{ MFG,2}$].
Let $\mu\in\cP_{m_0} 
(\Gamma^{\rm Lip}_{C,V})$ be a relaxed equilibrium and $(u,m)$ be the related mild solution. 
\begin{itemize}
\item[(a)] The function $u$ is locally semi-concave in $(J_i\setminus\{\vertex\})\times(0,T)$ in the sense of \cite[Definition 2.1.1]{MR2041617}. 
\item[(b)] Lemma~\ref{lemma:note4.9} holds replacing $\Gamma_t^{\textrm{opt}}[x]$ with $\Gamma_t^{\mu,\textrm{opt}}[x]$;
in particular, $D_x u(x,t)$ exists if and only if $\mathcal{A}(x,t)$ is a singleton and, in this case, equation~\eqref{eq:10003} holds with $\ell_i(\cdot,\cdot,t) $
replaced by~$L_i[m(t)]$.

\end{itemize}
\end{lemma}

In Proposition~\ref{prp:CCC2_4.1+4.2} below, we will see that $u$ is a bilateral subsolution (see \cite[Definition III.2.27]{BCD}) of the Hamilton-Jacobi equation and is differentiable at least at the points~$(x,t)$ such that $x$ belongs to the support of $m(t)$ and does not coincide with~$\vertex$. To this end, some new notation is useful. Set
\begin{eqnarray*}
\textrm{supp}(m(t))&=&\left\{x\in\cG\;:\; \forall \omega, \textrm{ open neighborhood of }x,\quad m(t)(\omega)>0\right\}\\
Q_m&=&\left\{(x,t)\in\cG\times(0,T)\;:\; x\in(\cG\setminus\{\vertex\})\cap \textrm{supp}(m(t))\right\}\\
\partial Q_m&=&\left\{(\vertex,t)\;:\; t\in(0,T ),\quad \vertex\in \textrm{supp}(m(t))\right\}
\end{eqnarray*}
and introduce the superdifferential of $u$ at $(x,t)\in Q_m$ as
\begin{equation*}
D^+_{x,t}u(x,t)=\left\{(\pi,q)\in \R^2\;:\; \limsup_{y\to x,\theta\to t}\frac{u(y,\theta)-u(x,t)-\pi(\theta-t)-q(\bar y-\bar x)}{|\bar y-\bar x|+|\theta-t|}\leq 0\right\},
\end{equation*}
where $x=\bar x e_j$, $y=\bar ye_j$ (note that $j$ is uniquely defined, from the definition of~$Q_m$). The set $D^+_{x,t}u(x,t)$ is nonempty, because, from Lemma \ref{lemma:reg_u_mfg}, $u$ is locally semi-concave in $(x,t)$, see~\cite[Proposition 3.3.4]{MR2041617}. 

\begin{proposition}\label{prp:CCC2_4.1+4.2}
Let us make the assumptions [$H_{ MFG,1}$] and [$H_{ MFG,2}$].
Let $\mu\in\cP_{m_0}(\G)$ be a relaxed equilibrium and $(u,m)$ be the associated mild solution. For any $(x,t)\in Q_m$, 
\begin{itemize}
\item[$(a)$] 
$
-p_1+H(x,t,p_2)=0$, for every $(p_1,p_2)\in D^+_{x,t}u(x,t)$
\item[$(b)$] $u$ is differentiable at $(x,t)$.
\end{itemize}
\end{proposition}
\begin{proof}
The proof follows the same arguments as the proof of \cite[Proposition 3.16]{AMMT3} ; therefore, we sketch them and refer the reader to  \cite{AMMT3} for the details.\\
$(a)$. Fix $(x,t)\in Q_m$ and consider $(p_1,p_2)\in D^+u(x,t)$. 
There exists $i, \; 1\le i \le N$,  such that  $x\in J_i$ with $x=\bar x e_i$. 
By Proposition~\ref{mild_sol_u}, $u$ is a viscosity solution to problem~\eqref{HJ} where the Hamiltonians $H_i$ and $H_\vertex$ are the ones given at the beginning of Section \ref{SHJB} replacing $\ell_i(\cdot,\cdot,t)$ and $\ell_0(t)$ respectively with $L_i[m(t)]$ and $L_0[m(t)]$. Since $u$ is a viscosity subsolution to~\eqref{HJ}, we have: $-p_1+H(x,t,p_2)\leq 0$. We want to prove the reverse inequality. Since $x\in\textrm{supp}(m(t))$, there exists a trajectory $\gamma\in \GammamuoptC[\gamma(0)]$ with $\gamma(t)=x$. Let $r$ be small enough such that $\gamma(t-s)\in J_i\setminus\{O\}$ for every $s\in[0,r]$. The definition of the subdifferential and Remark~\ref{basic_facts_value} entail
\begin{equation*}
\int_{t-s}^t L_i[m(\tau)](\gamma(\tau),\gamma'(\tau))d\tau\leq -p_1s-p_2\int_{t-s}^t\gamma'(\tau)d\tau+o(r).
\end{equation*}
The regularity of~$m$ (see Remark~\ref{rmk:reg_m}) and of~$\gamma$ in~$(t-s,t)$ (see Lemma~\ref{lemma:EL}, point 4) implies that $L_i[m(\tau)](\gamma(\tau),\gamma'(\tau)) = L_i[m(t)](x, \gamma'(t))+o(1)$ for any $\tau\in(t-s,t)$. This implies
\begin{equation*}
0\leq -p_1-p_2\gamma'(t)\cdot e_i -L_i[m(t)](x,\gamma'(t))\leq -p_1+H_i(x,p_2,t) 
\end{equation*}
where the last inequality comes from the definition of~$H_i$.\\
\noindent $(b)$. It suffices to follow the arguments in the proof of \cite[Proposition 4.2]{CCC2} replacing \cite[Theorem 4.1]{CCC2} and 
\cite[Corollary 4.1]{CCC2} respectively with  point~$(a)$ and Lemma~\ref{lemma:reg_u_mfg}-$(a)$.
\end{proof}

\begin{remark}\label{rem:nodirac}
Under assumptions [$H_{ MFG,1}$], [$H_{ MFG,2}$] and the counterpart of $[H6]$ in Remark~\ref{rmk:h7},  the properties listed in Remark~\ref{rmk:h7} hold, and, furthermore, point masses in $m$ cannot appear/vanish in the interior of a given edge, see \cite[Theorem 3.17]{AMMT3}, \cite[Corollary 3.18]{AMMT3} and \cite[Proposition 3.19]{AMMT3}.
\end{remark}

\subsection{The continuity equation}
\label{cont_eq}
Consider a mild solution $(u,m)$ associated with some relaxed  equilibrium~$\mu\in \cPmzeroGammaLipCV$. Our aim is to show that $m$ satisfies a continuity equation in a suitable weak sense, see point 3. in Theorem \ref{thm:cont_eq}. 
The main difficulty is to give a rigorous weak formulation of the transmission condition at the vertex $0$, accounting for the conservation of mass and involving fluxes associated with the different edges. The latter condition will be  sufficiently general to describe, for example, the evolution of a  Dirac mass arriving at the vertex at time $t_-$,  as well as the creation of a Dirac mass at the vertex. Remark \ref{rem_characterization} below contains a discussion on open problems related to this formulation.
\\
Hereafter, $\cD'(0,T)$ stands for the space of distributions in $(0,T)$.
\begin{theorem}\label{thm:cont_eq}
Let us make the assumptions [$H_{ MFG,1}$] and [$H_{ MFG,2}$].
Let $\mu\in \cPmzeroGammaLipCV$ be a relaxed equilibrium and $(u,m)$ be the related mild solution.
\begin{enumerate}
\item   For every $i=1,\dots, N$,  consider a function $\psi \in C^\infty(\cG)$, supported in $ J_i\setminus\{\vertex\}$, such that $\psi(x)=1$ for all $x\in J_i$ such that $d(x,\vertex)\geq 1$ and $\psi|_{J_i}$ is nondecreasing with respect to $d(x,\vertex)$.
For $\epsilon>0$, set $\psi_\epsilon(x)= \psi(\frac x \epsilon)$. The family of functions 
\begin{equation*}
%\begin{array}{rc}
t\mapsto -\int_{J_i} D_x\psi_\epsilon(x)\partial_pH(x,D_xu(x,t),t) m(t)(dx) 
\end{equation*}
converges  in $\cD'(0,T)$ as $\epsilon\to 0$ to some distribution $q_i$, which does not depend on the choice of $\psi$.
\item The following identities hold in  $ {\mathcal D}'(0,T)$ with $q_i$ obtained in point 1.:
\begin{eqnarray}
\label{eq:defqi}  
\;\;\frac d {dt} \Bigl( m(\cdot)\left(J_i\setminus\{\vertex\}\right) \Bigr) =\frac d {dt} \left(\int_\G \car_{\gamma(\cdot)\in J_i\setminus\{\vertex\}} \mu(d\gamma)\right) &=&
q_i.
\end{eqnarray}
\item  For every $\phi\in C^\infty(\cG\times[0,T])$ such that $\textrm{supp}(\phi(\cdot,t))$ is contained in a compact subset of~$\cG$ independent of~$t$, the following identities hold in $\cD'(0,T)$: for all $i=1,\dots,N$,
\begin{multline*}%\label{eq:continuity22}
\frac d {dt} \left(\int_\cG \car_{x\in J_i\setminus\{\vertex\}} \phi(x,\cdot) m(\cdot)(dx)\right)=\phi(\vertex,\cdot) q_i\\
+ \int_{\cG}\car_{x\in J_i\setminus\{\vertex\}}\left[\partial_t\phi(x,\cdot)-\partial_pH(x,D_xu(x,\cdot),\cdot)D\phi(x,\cdot)\right]m(\cdot)(dx),
\end{multline*}
and
\begin{eqnarray}
\label{eq:continuity}
&&\qquad \frac{d}{dt}\left[\int_{\cG}\phi(x,\cdot)m(\cdot)(dx)\right]=m(\cdot)(\{\vertex\})\partial_t\phi(\vertex,\cdot)\\ \notag &&
+\sum_{i=1}^N\int_{\cG}\car_{x\in J_i\setminus\{\vertex\}}\left[\partial_t\phi(x,\cdot)-\partial_pH(x,D_xu(x,\cdot),\cdot)D\phi(x,\cdot)\right]m(\cdot)(dx),\\
\label{eq:claim2_EC}
&&\qquad \frac d {dt}\left[m(\cdot)(\{\vertex\})\right]+\sum_{i=1}^N q_i=0.
\end{eqnarray}
\end{enumerate}
\end{theorem}
\begin{proof}
The proof follows the same arguments as that of~\cite[Theorem 3.20]{AMMT3}. Since it is delicate and long, we do not intend to reproduce it here. It mainly consists of observing that  \cite[Lemma 3.22]{AMMT3} (on the boundedness and measurability of the function $t\mapsto m(t)(\{x\})$ for all $x\in \cG$), and \cite[Lemma 3.23]{AMMT3} (on the boundedness and measurability of the function $t\mapsto  \ds \int_\cG \car_{x\in J_i\setminus\{\vertex\}} \phi(x,\cdot) m(\cdot)(dx)$) 
 hold in the present setting. The remainder of the proof is exactly the same as that of \cite[Theorem 3.20]{AMMT3}using, respectively, Proposition~\ref{prp:CCC2_4.1+4.2} and Lemma~\ref{lemma:reg_u_mfg}-(b) instead of \cite[Proposition 3.16-(b)]{AMMT3} and \cite[Lemma 3.14-(b)]{AMMT3}.
The latter results lead to the key observation that
$ \ds
\car_{\gamma(t)\ne 0}\gamma'(t) = - \car_{\gamma(t)\ne 0} \partial_p H \left(\gamma(t), D_xu(\gamma(t) ,t),t \right)$, 
for all $t\in [0,T]$ and for  $\mu-\textrm{a.e. } \gamma\in \Gamma$.
\end{proof}

\begin{remark}
Equation \eqref{eq:continuity}  in particular implies that
for all $\phi\in C_0^\infty (\cG\times (0,T) )$,
\begin{multline*}
\int_0^T m(t)(\{\vertex\})\partial_t\phi(\vertex,t) dt \\
+\sum_{i=1}^N   \int_0^T \int_{\cG}\car_{x\in J_i\setminus\{\vertex\}}\left[\partial_t\phi(x,t)-\partial_pH(x,D_xu(x,t),t)D_x\phi(x,t)\right]m(t)(dx) dt =0.
\end{multline*}
\end{remark}

\begin{remark}\label{rem_characterization}
Here, we do not deal with the uniqueness of mild solutions, that can be obtained under suitable assumptions, in particular separate dependence of the costs w.r.t. $\alpha$ and $m$, and monotonicity w.r.t. $m$, see \cite[Subsection 3.7]{AMMT3}. 

An open question is whether the mild solutions are characterized by
the system of PDEs coupling the Hamilton-Jacobi equation  and the continuity equation, respectively obtained in Proposition \ref{mild_sol_u} and Theorem \ref{thm:cont_eq}.
\end{remark}

\section{General networks}\label{sec_general_nets}
We now briefly discuss how to extend the results of the previous sections to a general network.
Hereafter, a {\it network} $\mathcal N$ is a connected subset of $\mathbb R^d$, made up of a finite number of non-self-intersecting closed lines, called {\sl edges},
which connect the nodes, called {\sl vertices}.  The finite collections of vertices and edges are respectively denoted by ${\mathcal V}=\{O_1,\dots,O_n\}$ and ${\mathcal E}=\{J_1,\dots, J_N\}$. We assume that for $i,j\in\{1,\dots,N\}$ with $i\ne j$, $J_i\cap J_j$ is empty or made of a single vertex. The length of the edge $J_i$ may be finite or infinite and is denoted by $l_i$ with $l_i\in(0,\infty]$. The edges can be arbitrarily oriented, but, for simplicity, we choose that an edge~$J_i$ connecting two vertices $O_{j_1}$ and $O_{j_2}$, with $j_1<j_2$, is oriented from $O_{j_1}$ toward $O_{j_2}$. This induces a natural parametrization $\pi_i:[0,l_i]\to J_i$ with $\pi_i(s)=[(l_i-s)O_{j_1}+ s O_{j_2}]/l_i$ for $s\in [0,l_i]$. Similarly, an infinite edge~$J_i$ starting from vertex~$O_j$ has the natural parameterization $\pi_i:[0,\infty)\to J_i$ with $\pi_i(s)=O_{j}+ s e_i$ for $s\in [0,\infty)$, where $e_i$ is the unit vector giving the direction and orientation of
$J_i$. Finally, for a given vertex $O\in{\mathcal V}$, we define the  {\sl junction related to $O$} as the union of $\{O\}$ and
all the relative interiors of the edges that are adjacent to $O$.

\subsection{A class of optimal control problems on the general network $\mathcal N$}
To any trajectory $y\in C([t, T]; \mathcal N)\cap W^{1,p} (t, T; \R^d)$ with $y(t)=x$, we associate the cost $J_t(x,y)$ as in~\eqref{eq:46}, with the running cost
\begin{eqnarray*}
L(y,\alpha,\tau)&=&\sum_{i=1}^N \ell_i(\pi_i^{-1}(y),\alpha\cdot e_i, \tau) \car_{y\in J_i\setminus{\mathcal V}}+ \sum_{j=1}^n\ell_{O_j}(\tau)\car_{y=O_j},\\
\ell_{O_j}(\tau)&=&\min\left\{\ell_{*j}(\tau),\min_{i=1,\dots,N:\ O_j\in J_i}\ell_i(\pi_i^{-1}(O_j),0,\tau)\right\}
\end{eqnarray*}
and the terminal cost
\begin{equation*}g(y)= \sum_{i=1}^N g_i(\pi_i^{-1}(y))\car_{y\in J_i\setminus{\mathcal V}}+\sum_{j=1}^n\min\left\{g_{*j},\min\limits_{i=1,\dots,N:\ O_j\in J_i}g_i(\pi_i^{-1}(O_j))\right\}\car_{y=O_j},
\end{equation*}
where $\ell_i$ fulfills the assumptions made in Section \ref{sec:OC},
while $\ell_{*j}$, $g_{*j}$ and $g_i$ fulfill, respectively, all the assumptions made in Section \ref{sec:OC} about $\ell_{*}$, $g_{*}$ and $g$.
With these assumptions, we can follow the same course as in Section \ref{sec:OC}. The following theorem summarizes the results that are obtained:
\begin{theorem}\label{thm:contr_N}
Under the above assumptions about running and terminal costs, each statement in Section~\ref{sec:OC} is valid in the present setting.
\end{theorem}

In the next two paragraphs, we briefly indicate how one can extend Proposition~\ref{prp:ex_OT} and Theorem~\ref{optcur_lip} respectively devoted to the existence and to the Lipschitz continuity of optimal trajectories. All other results can be extended in a straightforward manner.

\paragraph{Existence of optimal trajectories}
The statement in Proposition~\ref{prp:ex_OT} about the existence of optimal trajectories holds for the general networks defined above.

\begin{proposition}
Under assumptions [H1]-[H4] (with obvious adaptations), for each  $(x,t)\in \cG\times[0,T]$, there exists an optimal trajectory.
\end{proposition}

Let us sketch the proof of this result.
For $t=T$, there is nothing to do. Consider $(x,t)\in{\mathcal N}\times (0,T)$ and a sequence $(y_n)_{n\in{\mathbb N}}$
such that $y_n\in\Gamma_t[x]$ and $u(x,t)=\lim_{n\to\infty}J_t(x;y_n)$. Arguing as in the proof of Proposition~\ref{prp:ex_OT}, possibly after the extraction of a subsequence, there exists $y\in\Gamma_t[x]$
such that $y_n$ converges to~$y$ in $W^{1,p} (t, T; \R^d)$ weakly and in $C([t,T],{\mathcal N})$.
To prove that $y$ is an optimal trajectory, it is enough to prove inequality~\eqref{LSC2}.
Note that for $n$ sufficiently large, $y_n$  visits no other vertices than those visited by $y$, from the uniform convergence of $y_n$ to $y$.
We argue differently according to the number of vertices visited by $y$.\\
If in $[t,T]$, $y$ visits one vertex at most, then, for $n$ sufficiently large,
each $y_n$ visits at most one vertex (in fact none if $y$ does not visit any vertex, and
at most one vertex that is the one visited by $y$ if $y$ does visit a vertex). In this case, we can argue exactly as in the proof of Proposition~\ref{prp:ex_OT}.\\
Assume now that $y$ visits two or more different vertices. In this case, $\delta=\min_{i=1,\dots,N}{l_i}$ is finite. 
As mentioned above, for $n$ sufficiently large, each $y_n$ does not visit any vertex that is not visited by $y$.
Without loss of generality, we can relabel the vertices as follows: let us set \begin{equation*}
\tau_1=\mathrm{argmin}\{\tau\in[t,T]:\; y(\tau)\in{\mathcal V}\},\quad\textrm{ and  }y(\tau_1)=O_1.
\end{equation*}
 Because $y$ is H\"older continuous and $\delta$ is strictly positive, there exists a finite integer $k$ such that  for $j=2,\dots, k$, we can define
\begin{equation*}
\tau_j=\mathrm{argmin}\{\tau\in[\tau_{j-1},T]:\; y(\tau)\in{\mathcal V}\setminus\{O_{j-1}\}\},\quad\textrm{and }y(\tau_j)=O_j,
\end{equation*}
and $y|_{ [\tau_k, T]}$ does not cross any vertex other than $O_k$.
Hence, $y|_{(\tau_{j-1},\tau_j)}$ cannot visit any other vertex than
$O_{j-1}$ (it may  visit  $O_{j-1}$ even infinitely many times).

For $j=2,\dots, k$, let us introduce
\begin{equation*}
\bar \tau_{j-1}=\mathrm{argmax}\{\tau\in[\tau_{j-1},\tau_{j}]:\; d(y(\tau), O_j)=\delta/2\},
\end{equation*}
and, for simplicity of notation, set $\bar \tau_0=t$ and $\bar \tau_{k}=T$. Note that for $j\in \{2,\dots, k\}$, $d(y(\tau), O_j)<\delta/2$ in $(\bar \tau_{j-1},\tau_j]$.
Clearly, for any $j=1,\dots, k$, $y|_{[\bar \tau_{j-1},\bar \tau_j)}$ takes its values in the junction related to $O_j$,
and so does $y_n|_{(\bar \tau_{j-1},\bar \tau_j)}$ for $n$ sufficiently large. \\
In order to prove the inequality~\eqref{LSC2}, we first observe that
\begin{equation*}
\liminf_{n\to\infty} \int_t^T  L(y_n(\tau),y'_n(\tau),\tau) d\tau \geq
\sum_{j=1}^{k}\liminf_{n\to\infty} \int_{\bar \tau_{j-1}}^{\bar \tau_{j}}L(y_n(\tau),y'_n(\tau),\tau) d\tau.
\end{equation*}
We may then apply the Lemma~\ref{lem:LSC} to each term on the right-hand side, replacing the interval $[t,T]$ with $[\bar \tau_{j-1},\bar \tau_{j}]$. We obtain 
\begin{equation*}
\liminf_{n\to\infty} \int_{\bar \tau_{j-1}}^{\bar \tau_{j}}L(y_n(\tau),y'_n(\tau),\tau) d\tau\geq \int_{\bar \tau_{j-1}}^{\bar \tau_{j}}L(y(\tau),y'(\tau),\tau) d\tau.
\end{equation*}
We conclude the proof of~\eqref{LSC2} by injecting the latter inequalities into the former.

\paragraph{ Local Lipschitz continuity of the optimal trajectories}
The statement in Theorem~\ref{optcur_lip} on the local Lipschitz continuity of optimal trajectories holds for the general networks defined above.

\begin{proposition}
Under assumptions [H1] to [H5] (with obvious adaptations), for every $(x,t)\in \cG\times[0,T]$ and for every $y\in \Gamma^{\rm{opt}}_t[x]$, $y'$ belongs to $L^\infty(t,T;\R^d)$. Moreover, there exists a positive constant $V$ (depending on $x$ and $T-t$, and that may blow up as $t\to T$) such that
$\|y'\|_\infty \leq V$.
\end{proposition}

Let us sketch the proof.
Consider $(x,t)\in{\mathcal N}\times[0,T)$ and $y\in \Gamma_t^{\textrm{opt}}[x]$. 
We are going to apply the same arguments as in the proof of Theorem~\ref{optcur_lip} edge by edge, proceeding from the last to the first edge visited by $y$.\\
Introduce the sets $I_{O_j}=\{\tau\in[t,T]:\; y(\tau)=O_j\}$ for $j=1,\dots,n$ and $I_{J_i}=\{\tau\in[t,T]:\; y(\tau)\in J_i\setminus{\mathcal V}\}$ for $i=1,\dots,N$.
For the sets $I_{O_i}$, we can just invoke Stampacchia's Theorem as in point $(1)$ in the proof of Theorem~\ref{optcur_lip}.

Let us distinguish two cases, whether $y(T)$ coincides or not with a vertex.
\begin{description}
\item{\it Case $1$: $y(T)\notin {\mathcal V}$.}
  There exist $t_1\in[t,T)$ and $i\in \{1,\dots,N\}$ such that $(t_1,T]\subset I_{J_i}$ and either $y(t_1)\in {\mathcal V}$ or $t_1=t$. Following the arguments in
point $(2)$ in the proof of Theorem~\ref{optcur_lip}, we obtain the desired result in the interval $[t_1,T]$.\\
We are left with obtaining the estimate in the interval $[t,t_1]$. If there exists $j\in\{1,\dots, N\}$ such that $[t,t_1)\subset I_{J_j}$, then the same arguments as those in point $(4)$ of the proof of Theorem~\ref{optcur_lip} lead to the desired result.\\
Suppose now that for all  $j\in \{1,\dots,N\}$, $[t,t_1)\not \subset I_{J_j}$. If, in a left neighbourhood of $t_1$, the trajectory $y$ visits the vertex $y(t_1)$ again (possibly infinitely many times), then for each interval $[\overline t,\underline t]$ with $y(\overline t)=y(\underline t)=y(t_1)$  the arguments in point~$(3)$ in the proof of Theorem~\ref{optcur_lip} yield the result. Hence, we can focus on the case where, for some $t_2\in[t,t_1]$, there holds: $y(t_2)\in{\mathcal V}$, $y(t_2)\ne y(t_1)$ and $(t_2,t_1)\subset I_{J_j}$. In this case, the trajectory crosses the entire edge $J_j$ in the interval $[t_2,t_1]$. Fix any point $\tilde x\in J_j\setminus{\mathcal V}$. There exists $\tilde t_2\in(t_2,t_1)$ such that $y(\tilde t_2)=\tilde x$. With the same arguments as in point $(4)$ in the proof of Theorem~\ref{optcur_lip}, we obtain the estimate in the interval $[\tilde t_2,t_1]$. Then, using this information, arguments similar to those in point~$(2)$ produce the estimate in the interval $[t_2,\tilde t_2]$. Consequently, we have obtained the desired estimate in $[t_2,t_1]$. Note that, since $y$ is H\"older continuous 
from Proposition \ref{prp:ex_OT}, $t_1-t_2$ is bounded away from zero. Repeating this argument a finite number of times, we obtain the desired estimate in the interval $[t,t_2]$. This concludes the proof.
\item{Case $2$: $y(T)\in{\mathcal V}$.} This case is simpler than the latter, and only some of the arguments are used. We omit the details.
\end{description}

\subsection{Mean Field Games on the general network $\mathcal N$}
Let $m_0$ be a fixed probability measure on $\mathcal N$ with compact support.
For $i=1,\dots, N$, let us define the running costs~$L_i$ and the terminal costs~$G_i$ as in Section~\ref{subsec:setting}, replacing $\cG$ by $\mathcal N$. Similarly, for each $j=1,\dots,n$, the  functions $L_{*j}$ and $G_{*j}$ are defined respectively as the functions $L_{*}$ and $G_{*}$ of Section~\ref{subsec:setting}, replacing  $\cG$ with $\mathcal N$.\\
To any trajectory $y\in C([t, T]; \mathcal N)\cap W^{1,p} (t, T; \R^d)$ with $y(t)=x$, we associate the cost $J_t(x,y)$ as in~\eqref{eq:46} with 
\begin{eqnarray*}
L[m](y(t),y'(t))&=&\sum_{i=1}^N L_i[m](\pi_i^{-1}(y(t)),y'(t)\cdot e_i) \car_{y(t)\in J_i\setminus{\mathcal V}}+ \sum_{j=1}^n L_{O_j}[m]\car_{y(t)=O_j},\\
L_{O_j}[m]&=&\min\left\{L_{*j}[m],\min_{i=1,\dots,N:\ O_j\in J_i}L_i[m](\pi_i^{-1}(O_j),0)\right\}
\end{eqnarray*}
and 
\begin{equation*}G[m](x)= \sum_{i=1}^N G_i[m](\pi_i^{-1}(x))\car_{x\in J_i\setminus{\mathcal V}}+\sum_{j=1}^n\min\left\{G_{*j}[m],\min\limits_{i=1,\dots,N:\ O_j\in J_i}G_i[m](\pi_i^{-1}(O_j))\right\}\car_{x=O_j}.
\end{equation*}
The assumptions made here on the costs are straightforward adaptations of [$H_{MFG,1}$] and 
[$H_{MFG,2}$]. They allow us to  follow the same course as in Section~\ref{sect:MFGequil} and obtain similar results:
\begin{theorem}
Under the assumptions above, each statement in Section~\ref{sect:MFGequil} holds in the present setting, (with natural adaptations).
\end{theorem}
\begin{proof}
We skip the proof, since the arguments are similar to those given in Section~\ref{sect:MFGequil}.
\end{proof}

\appendix

\section{Proof of Lemma \ref{lemma:barL0} } %and \ref{lemma:prp_cal_J}}
\label{sec:app1}

\begin{proof}[Proof of  Lemma \ref{lemma:barL0}]
  
$i)$. By the definition of the convex envelope, it is enough to prove that $\mathcal L_0(0,s)\geq \ell_\vertex(s)$ for every $s\in[0,T]$. To this end, we first note that any vector $\alpha \in\R^N$ can be expressed as a convex combination of vectors in $\cup_{i=1}^N\R e_i$. Therefore,  $\mathcal L_0(\alpha,s)$ is finite for any $(\alpha,s)\in \R^N\times [0,T]$, i.e. the domain of $\mathcal L_0$ is $ \R^N\times [0,T]$.\\
From \cite[Theorem 2.1]{MR986024}, there exist $\alpha_1,\dots,\alpha_{N+1}\in \R^N$ and $\lambda_1,\dots,\lambda_{N+1}\in (0,1]$ such that $\mathcal L_0(0,s)=\sum_{i=1}^{N+1}\lambda_i \min\left(\ell(\vertex,\alpha_i,s),\ell_\vertex(s)+\chi_{0}(\alpha_i)\right)$,   $\sum_{i=1}^{N+1}\lambda_i=1$ and $\sum_{i=1}^{N+1}\lambda_i\alpha_i=0$.
Possibly after a permutation of the indices, we may suppose that there exists $\bar N\in\{0,\dots,N+1\}$ such that
\begin{equation*}
\min\left(\ell(\vertex,\alpha_i,s),\ell_\vertex(s)+\chi_{0}(\alpha_i)\right)=
\left\{\begin{array}{ll}
\ell(\vertex,\alpha_i,s)&\qquad \textrm{if }  1 \le i \le \bar N, 
\\
\ell_\vertex(s)+\chi_{0}(\alpha_i)&\qquad \textrm{if } \bar N<  i \le N+1.
\end{array}\right.
\end{equation*}
From the definition of~$\chi_{0}$, we infer that $\alpha_i=0$ if $ \bar N<i \le N+1$. This and \eqref{eq:460} imply that $\min\left(\ell(\vertex,s,\alpha_i),\ell_\vertex(s)+\chi_{0}(\alpha_i)\right)=\ell_\vertex(s)$ if $\bar N < i\le  N+1$.
We now distinguish three cases.
\begin{enumerate}
\item $\bar N=0$: in this case, $\alpha_i=0$ for $i=1,\dots,N+1$ and $\mathcal L_0(0,s)=\sum_{i=1}^{N+1}\lambda_i \ell_\vertex(s)=\ell_\vertex(s)$.
\item  $\bar N=N+1$: since $\ell_\vertex(s)\le \min_{i} \ell_i(0,0,s)$, we may assume without loss of generality that $\alpha_i\ne 0$  for all $i=1,\dots, N+1$. Hence, $\mathcal L_0(0,s)=\sum_{i=1}^{N+1}\lambda_i \ell(\vertex,\alpha_i,s)$, and $\sum_{i=1}^{N+1}\lambda_i\alpha_i=0$. Since $\mathcal L_0(0,s)$ is finite, we deduce from the definition of $\ell$ that for all $i\in \{1,\dots, N+1\}$, $\alpha_i \in \cup_{j=1}^N \R e_j$. Possibly after another permutation of the indices, we may assume that there exist $\bar N_1,\dots,\bar N_N\in\{0,\dots,N+1\}$ such that $\bar N_1\leq \bar N_2\leq \dots\leq \bar N_N$ and $\alpha_i\in \R e_j$ for $\bar N_j\le i< \bar N_{j+1}$. Hence, $ \mathcal L_0(0,s) =\sum_{j=1}^N\sum_{i=\bar N_j}^{\bar N_{j+1}-1}\lambda_i \ell_j(0,\alpha_i\cdot e_j,s)$,
$\sum_{j=1}^N\sum_{i=\bar N_j}^{\bar N_{j+1}-1} \lambda_i=1$ and   $\sum_{j=1}^N\sum_{i=\bar N_j}^{\bar N_{j+1}-1} \lambda_i\alpha_i=0$. But, since the vectors $e_j$ are linearly independent, the latter identity implies that for all $j\in \{1,\dots, N\}$,
$\sum_{i=\bar N_j}^{\bar N_{j+1}-1} \lambda_i\alpha_i =0 $. Set $\bar \lambda_j= \sum_{i=\bar N_j} ^{\bar N_{j+1} -1}\lambda_i$.
If  $\bar \lambda_j>0$,  the convexity of $\ell_j$ yields the following:
\begin{equation*}
\mathcal L_0(0,s)=
\sum_{j=1}^N\bar \lambda_j \sum_{i=\bar N_j}^{\bar N_{j+1}-1}\frac{\lambda_i}{\bar \lambda_j} \ell_j(0,\alpha_i\cdot e_j,s)\geq \sum_{j=1}^N\bar \lambda_j\ell_j(0,0,s).
\end{equation*}
Since $\sum_{j=1}^N\bar \lambda_j=1$, this implies that
$\cL_0(0,s)\geq \sum_{j=1}^N\bar \lambda_j\ell_\vertex(s) = \ell_\vertex(s)$.
\item $ 0< \bar N\le N $:  In this case, $\alpha_i=0$ for $\bar N < i\le N+1$, and as in the previous case,  we may assume without loss of generality that $\alpha_i\ne 0$  for all $i=1 ,\dots,  \bar N$. Setting $\bar \lambda_0= \sum_{i=\bar N +1} ^{N+1} \lambda_i$, we see that $\frac {\mathcal L_0(0,s) -   \bar \lambda_0 \ell_\vertex(s)} {1-\bar \lambda_0} =    \sum_{i=1}^{\bar N }\frac {\lambda_i}{1-\bar \lambda_0} \ell(\vertex,\alpha_i,s)$, and $\sum_{i=1}^{\bar N} \frac {\lambda_i}{1-\bar \lambda_0}  \alpha_i=0$.
Arguing as in the previous case yields $  \sum_{i=1}^{\bar N }\frac {\lambda_i}{1-\bar \lambda_0} \ell(\vertex,\alpha_i,s) \ge
\ell_\vertex(s)$. Combining the observations above, we find that $\mathcal L_0(0,s)\ge  \ell_\vertex(s)$.
\end{enumerate}

\noindent
$ii)$  This part of the proof borrows some ideas from \cite[Theorem 2.1]{MR986024}.  Fix $(\alpha,s)\in \R^N\times [0,T]$. We aim at proving that
\begin{equation}\label{eq:liminf_L0}
{\mathcal L}_0(\alpha,s)\leq \liminf_{\tau\to s}{\mathcal L}_0(\alpha,\tau).
\end{equation}
Observe that the functions $(\alpha, \tau) \mapsto \ell (\vertex,\alpha,\tau)$ and $(\alpha,\tau)\mapsto \ell_\vertex(\tau)+\chi_{0}(\alpha)$ are LSC in $\R^N\times [0,T]$. Let us denote by $D_\ell$ and $D_0$ the domains of $\ell(\vertex,\cdot,\cdot)$ and of $(\alpha,\tau)\mapsto  \ell_\vertex(\tau)+\chi_{0}(\alpha)$ respectively. We know that $D_{\ell}=(\cup_{i=1}^N \R e_i )\times [0,T]$ and $D_0 =\{0\}\times [0,T]$, whereas ${\rm{dom}} (\cL_0)= \R^N\times [0,T]$.
\\
Let the sequence $(\tau_n)_n$ be such that $\lim_{n\to\infty}  {\mathcal L}_0(\alpha,\tau_n)= \liminf_{\tau\to s}{\mathcal L}_0(\alpha,\tau)$. From \cite[Theorem 2.1]{MR986024}, for every $n\in\N$, there exist $q_n\in\{1,\dots,N+1\}$, $\alpha^n_i\in D_l\cup D_0$ and $\lambda^n_i\in [0,1]$ for $i=1,\dots,q_n$, with $\sum_{i=1}^{q_n}\lambda^n_i=1$ such that
\begin{equation*}
{\mathcal L}_0 (\alpha,\tau_n)=\sum_{i=1}^{q_n}\lambda^n_i \min\left(\ell(\vertex,\alpha^n_i,\tau_n),\ell_\vertex(\tau_n)+\chi_{0}(\alpha^n_i)\right)\quad\textrm{and}\quad \alpha=\sum_{i=1}^{q_n}\lambda^n_i \alpha^n_i.
\end{equation*}
Note that $\min\left(\ell(\vertex,0,\tau_n),\ell_\vertex(\tau_n)+\chi_{0}(0)\right)=\ell_\vertex(\tau_n)$. As in the proof of \cite[Theorem 2.1]{MR986024}, possibly after the extraction of a subsequence, we may assume that $q_n=\bar q$ is an integer independent of $n$. Hence, ${\mathcal L}_0 (\alpha,\tau_n)=\sum_{i=1}^{\bar q}\lambda^n_i \min\left(\ell(\vertex,\alpha^n_i,\tau_n),\ell_\vertex(\tau_n)+\chi_{0}(\alpha^n_i)\right)$ and $ \alpha=\sum_{i=1}^{\bar q}\lambda^n_i \alpha^n_i$. Note that $(\alpha^n_i, \tau_n)\in  D_\ell\cup D_0$, which is a closed set.

Thus, after further extractions of subsequences, for each $i=1,\dots, \bar q$, there are three possible cases:
\begin{enumerate}
\item  there exists $\alpha_i\ne 0$ such that $\lim_n \alpha^n_i=\alpha_i$ and $(\alpha_i,s)\in D_\ell$
\item $\lim_n \alpha^n_i=0$
\item  $\lim_n |\alpha^n_i|=\infty$.
%\item  there exists $\alpha_i$ such that $\lim_n \alpha^n_i=\alpha_i$ and $(alpha_i,s)\notin D_\ell\cup D_0$.
\end{enumerate}
Hence, possibly after exchanging the indices, we may assume that there exist $q_\ell$ and $q$,  $0\le q_\ell\le q\le \bar q$, such that
\begin{eqnarray}
\lim_n \alpha^n_i=\alpha_i\in\bigcup_{i=1}^N\R e_i\setminus\{0\}\qquad &\textrm{if }& 1\leq i\leq q_\ell \label{2.9}\\
\lim_n \alpha^n_i=\alpha_i=0 \qquad&\textrm{if }& q_\ell+1 \leq i \leq q\label{2.9bis}\\
\lim_n |\alpha^n_i|=\infty \qquad&\textrm{if }& q+1 \leq i \leq \bar q.\label{2.11}
\end{eqnarray}
Then, possibly after a further extraction of subsequences, we may assume that there exists a map $j:\{1,\dots, \bar q\}\to \{1,\dots, N\}$ such that  if $1\le i\le q_\ell$ or if  $q+1\le i\le \bar q$, then 
$\alpha_i^n \in \R e_{j(i)}\setminus \{0\}$,
and  if  $q_\ell+1\le i\le q$, then $\alpha_i^n \in \R e_{j(i)}$. In particular, if  $q+1\le i\le \bar q$,
$\min \left(\ell(\vertex,\alpha_i^n, \tau_n), \ell_\vertex(\tau_n)+\chi_0(\alpha_i^n)\right)=
\ell(\vertex,\alpha_i^n, \tau_n)=\ell_{j(i)} (\vertex,\alpha_i^n\cdot e_j,\tau_n)$

There exist $(\lambda_i)_{1\le i\le \bar q}$ such that, after a further extraction, $\lim \lambda_i^n =\lambda_i$. Then $\lambda_i\in [0,1]$ and $\sum_{i=1}^{\bar q} \lambda_i=1$. The lower semi-continuity of $\ell$ implies that for all $i\in \{1,\dots, q_\ell\}$, 
\begin{displaymath}
\liminf_{n\to\infty}  \min\left(\ell(\vertex,\alpha^n_i,\tau_n),\ell_\vertex(\tau_n)+\chi_{0}(\alpha^n_i)\right)= \liminf_{n\to\infty}  \ell(\vertex,\alpha_i^n, \tau_n)\ge \ell(\vertex,\alpha_i, s).
\end{displaymath}
Similarly, the lower semi-continuity of $\ell$ and \eqref{eq:460} imply that, for all $i\in \{q_\ell+1,\dots, q\}$,
\begin{displaymath}\liminf_{n\to\infty}  \min\left(\ell(\vertex,\alpha^n_i,\tau_n),\ell_\vertex(\tau_n)+\chi_{0}(\alpha^n_i)\right)\ge 
\min  \left(   \ell(\vertex,0, s), \ell_\vertex(\tau)  \right)= \ell_\vertex(\tau).
\end{displaymath}
Finally, the growth assumption in [H1] yields 
$\min\left(\ell(\vertex,\alpha^n_i,\tau_n),\ell_\vertex(\tau_n)+\chi_{0}(\alpha^n_i)\right)=\ell(\vertex,\alpha^n_i,\tau_n)\to\infty$,  for $ q+1\leq i\leq \bar q$. This implies that for $n$ large enough,
\begin{equation} \label{2.15}
{\mathcal L}_0(\alpha,\tau_n)   \ge  \sum_{i=1}^{q}\lambda^n_i
\min\left(\ell(\vertex,\alpha^n_i,\tau_n),\ell_\vertex(\tau_n)+\chi_{0}(\alpha^n_i)\right).
\end{equation}
We also deduce from  assumption [H1] that 
\begin{equation}\label{2.16}
\sum_{i=q+1}^{\bar q}\lambda^n_i \ell(\vertex,\alpha^n_i,\tau_n)\leq {\mathcal L}_0(\alpha,\tau_n)+C_0.
\end{equation}
Now, since $\alpha=\sum_{i=1}^{\bar q} \lambda_i^1 \alpha_i^1$ with $\alpha_i^1\in D_\ell\cup D_0$ for all $i=1,\dots,\bar q$,
there exists a constant $\tilde C$ independent of $n$ such that
\begin{equation}
\label{2.18ter}
{\mathcal L}_0 (\alpha,\tau_n)\leq \sum_{i=1}^{\bar q}\lambda^1_i \min\left(\ell(\vertex,\alpha^1_i,\tau_n),\ell_\vertex(\tau_n)+\chi_{0}(\alpha^1_i)\right)\leq \tilde C.
\end{equation}
Plugging \eqref{2.18ter} in \eqref{2.16}, we easily deduce that $\lambda_i=0$ for $q+1\le i \le \bar q$ and that $\sum_{i=1}^q \lambda_i=1$.
From the superlinear growth of the costs $\ell_i$,
we also deduce $\lim_{n\to \infty} \sum_{i=q+1}^{\bar q} \lambda_i^n |\alpha^n_i|=0$. Hence $\alpha= \sum_{i=1}^q \lambda_i \alpha_i$.  Moreover, $\alpha_i=0$ if $i> q_\ell$. Taking the $\liminf$ in \eqref{2.15} yields $ \liminf_{n\to \infty}  {\mathcal L}_0(\alpha,\tau_n)   \ge  \sum_{i=1}^{q_\ell }\lambda_i
\ell(\vertex,\alpha_i,s)  +    \sum_{i=q_\ell+1}^{q }\lambda_i\ell_\vertex(s)$.
This implies \eqref{eq:liminf_L0}, since $\alpha$ is a convex combination of the vectors $(\alpha_i)_{i=1,\dots, q} $ with $\alpha_i=0$ if $i> q_\ell$. The proof of $ii)$ is complete.

$iii)$ The function $(\alpha, s) \mapsto \cL_0(\alpha, s) $ is locally bounded, lower semi-continuous with respect to $s$ (for $\alpha$ fixed) and convex with respect to $\alpha$ (for $s$ fixed), therefore continuous with respect to $\alpha$ (for $s$ fixed). Point $iii)$ is then a consequence of ~\cite[Theorem 8]{ursell1939some}.
\end{proof}

\section*{Acknowledgments} The first author is partially supported by the chair Finance and Sustainable Development and FiME Lab (Institut Europlace de Finance).
The first and third authors are  partially supported by ANR (Agence Nationale de la Recherche) through project COSS, ANR-22-CE40-0010-01. 
The second author is partially supported by INDAM-GNAMPA project codice CUP E53C25002010001, by the project PRIN, codice MUR 2022W58BJ5 ``PDEs and optimal control methods in mean field games, population dynamics and multi-agent models'' and by the KAUST Research Funding ``Free Boundary Problems, Mean Field Games, Crowd Motion and Lipschitz Learning: The Infinity-Laplacian in Action'' under award no. CRG2024-6430.6.

{\small
\bibliographystyle{siamplain}
\bibliography{references}
}
\end{document}